\documentclass[12pt]{article}
\usepackage{amssymb}
\usepackage{graphics}

\newcounter{cclaim}[section]
\renewcommand{\thecclaim}{\arabic{cclaim}}
\newenvironment{cclaim}{\refstepcounter{cclaim}
	\par\medskip\par\noindent{\it Claim~\thecclaim.}\rm}
{\par\medskip\par}

\newenvironment{subproof}{\par\noindent{\it Proof of Claim.}}
{$\Box$\par\smallbreak}

\newcommand{\hide}[1]{}

\newcommand{\N}{{\mathbb N}}

\newcommand{\R}{{\mathbb R}}

\newcommand{\Z}{{\mathbb Z}}

\newcommand{\LL}{{\cal L}}

\newcommand{\beq}{\begin{equation}}
\newcommand{\ee}{\end{equation}}

\renewcommand{\d}{\partial}

\newcommand{\al}{\alpha}
\newcommand{\be}{\beta}
\newcommand{\ga}{\gamma}
\newcommand{\de}{\delta}
\newcommand{\eps}{\varepsilon}
\newcommand{\vphi}{\varphi}
\newcommand{\la}{\lambda}
\newcommand{\om}{\omega}
\newcommand{\reff}[1]{(\ref{#1})}

\newcommand{\diag}{\mathop{\rm diag}\nolimits}

\newtheorem{thm}{Theorem}[section]

\title{
Bounded and Almost Periodic Solvability of  
 Nonautonomous Quasilinear Hyperbolic Systems
}

\newcounter{thesame}
\author{
 Irina Kmit
\ \ \ Lutz Recke \ \ \ Viktor Tkachenko\\ [3mm]
{\small
Institute of Mathematics, Humboldt University of Berlin, Germany}
\\
{\small and Institute for Applied Problems of Mechanics and Mathematics,
}
\\
{\small
National Academy of Sciences of Ukraine, Lviv, Ukraine}\\
{\small   E-mail:
{\tt kmit@mathematik.hu-berlin.de}}\\ [3mm]
{\small
Institute of Mathematics, Humboldt University of Berlin, Germany}
\\
{\small   E-mail:
{\tt recke@mathematik.hu-berlin.de}}\\[3mm]
{\small
Institute of Mathematics, National Academy of Sciences of Ukraine, Kyiv, Ukraine
}
\\
{\small   E-mail:
{\tt vitk@imath.kiev.ua}}
}
\date{}
\begin{document}
\maketitle

\begin{abstract}
The paper concerns boundary value problems
for general nonautonomous  first order quasilinear hyperbolic systems in a strip.
We construct   small  global   classical solutions, assuming
 that the right hand sides are small.  
In the case that  all data of the quasilinear problem are almost periodic,  we prove that
the bounded solution is also almost periodic.  
 For the nonhomogeneous version of a linearized problem,  we provide  stable dissipativity conditions
ensuring a unique bounded  continuous solution for any smooth right-hand sides. 
In the autonomous case, this solution is two times continuously differentiable.
In the nonautonomous case,  the continuous solution  is  differentiable under
additional dissipativity conditions, which are essential. 
A crucial ingredient of our approach is a perturbation theorem for  general linear 
hyperbolic systems. One of the  technical complications  we overcome
is the ``loss of smoothness'' property of hyperbolic PDEs. 
\end{abstract}

{\bf Key words:} nonautonomous  quasilinear  hyperbolic systems,
 boundary value problems, bounded classical solutions, almost periodic solutions, dissipativity conditions,
 perturbation theorem for linear problems

\section{Introduction}\label{sec:intr}
\renewcommand{\theequation}{{\thesection}.\arabic{equation}}
\setcounter{equation}{0}

\renewcommand{\theequation}{{\thesection}.\arabic{equation}}
\setcounter{equation}{0}

\subsection{Problem setting}\label{sec:setting}
We consider a
 first order quasilinear hyperbolic system
\begin{equation}\label{eq:1}
\partial_t V + A(x,t,V)\partial_x V + B(x,t, V)V = f(x,t),\quad  x\in(0,1), t\in \R,
\end{equation}
where  $V=(V_1,\ldots,V_n)$ and $f = (f_1,\dots,f_n)$ are  vectors of real-valued functions, and
$A=(A_{jk})$ and $B = (B_{jk})$ are $n\times n$-matrices 
 of real-valued functions. 
The matrix  $A$ is supposed to have  $n$ real eigenvalues $A_j(x,t,V)$ in a neighborhood of $V=0$
in $\R^n$ such that
\begin{eqnarray*}
A_1(x,t,V)>\dots>A_m(x,t,V)> 0>A_{m+1}(x,t,V)>\dots>A_n(x,t,V)
\end{eqnarray*}
for some integer $0\le m\le n$.  These assumptions imply that there exists
 a smooth and  nondegenerate $n\times n$-matrix  $Q(x,t,V) = \left(Q_{jk}(x,t,V)\right)$ 
such that  
\begin{eqnarray*}
 Q^{-1}(x,t,V) A(x, t,V)Q(x,t,V) = \diag(A_1(x,t,V),\dots,A_n(x,t,V)). 
\end{eqnarray*}

 We supplement the system (\ref{eq:1})  with the boundary conditions 
\beq\label{eq:2*}
\begin{array}{ll}
  U_{j}(0,t)= (RZ)_j(t) + h_j(t), \;\;\; 1\le j\le m,\; t\in \R,  \\ [1mm]
 U_{j}(1,t)= (RZ)_j(t)  + h_j(t), \;\;\; m< j\le n,\; t\in \R,
\end{array}
\ee
where $R$ is a 
(time-dependent)
 bounded linear operator,
$$
Z(t)=\left(U_1(1,t),\dots,U_{m}(1,t),U_{m+1}(0,t),\dots,U_{n}(0,t)\right),
$$
and
\beq\label{z*}
U(x,t) =  Q^{-1}(x,t,V)V(x,t).
\ee

The purpose of the paper is to establish conditions on the coefficients 
$A$, $B$, $f$, and $h$
 and
the boundary operator $R$ ensuring that the problem \reff{eq:1}--\reff{z*} has 
 a unique
 small  global   classical solution, which is  two times continuously differentiable.
 If the data in \reff{eq:1} and \reff{eq:2*} are almost periodic 
 (respectively, periodic) in $t$,
 we  prove that the bounded solution is
almost periodic (respectively, periodic) in $t$ also.

Let
$$
\Pi = \{(x,t)\in\R^2\,:\,0 \le x \le 1\}
$$
and $BC(\Pi;\R^n)$
be the Banach space of
all continuous and bounded maps    
$u:\Pi \to \R^n$
with the usual $\sup$-norm
$$
\|u\|_{BC}=\sup\left\{|u_j(x,t)|\,:\,(x,t)\in\Pi, j\le n\right\}.
$$
Moreover,  $BC^k(\Pi;\R^n)$ denotes the space of $k$-times
continuously differentiable and bounded maps    
$u:\Pi \to \R^n$,
with  norm
$$\|u\|_{BC^k} =
\sum_{0\le i+j\le k}\left\|\partial_x^i\partial_t^ju\right\|_{BC}.
$$
We  also use the spaces $BC^k_t(\Pi;\R^n)$   of functions $u \in BC(\Pi;\R^n)$  such that
$\partial_t u,\dots,\partial^k_t u \in BC(\Pi;\R^n)$,
with  norm
$$
\|u\|_{BC^k_t} = \sum\limits_{j=0}^k\|\partial_t^j u\|_{BC}.
$$
Similarly, $BC^k(\R;\R^n)$ denotes  the space of  $k$-times
continuously differentiable and bounded maps $u:\R \to \R^n$. If $n=1$, we will simply write $BC^k(\R)$
for $BC^k(\R;\R)$, and likewise for all the spaces introduced above.

Given two Banach spaces $X$ and $Y$, the space of all  bounded linear operators $A: X \to Y$ is
denoted by $\mathcal{L}(X,Y)$, with the operator norm
$\|A\|_{\mathcal{L}(X,Y)}=\sup\{ \|Au\|_{Y}: u\in X, \|u\|_{X} \le 1\}.$
We will use also the usual notation $\LL(X)$ for $\LL(X,X)$.

 Let $\|\cdot\|$ denote the  norm in $\R^n$ defined by
$
\|y\|=\max\limits_{j\le n}|y_j|.
$
We suppose that the  data of the problem  (\ref{eq:1})--(\ref{z*})
satisfy the following conditions.

\begin{itemize}
  \item[\bf(A1)]
 There exists $\de_0 > 0$  such that
\begin{itemize}
\item
the entries of the matrices $A(x,t,V)$, $B(x,t,V)$, and  $Q(x,t,V)$ have bounded and
continuous partial derivatives  up to the second order in $(x,t)\in\Pi$ and in $V\in\R^n$
with $\|V\| \le \de_0$,
\item
there exists $\Lambda_0>0$ such that
$$
\begin{array}{ll}
\inf\left\{A_j(x,t,V)\,:\,(x,t)\in \Pi,\, \|V\| \le \de_0, 1\le j\le m\right\}\ge \Lambda_0,\\ [1mm]
\sup\left\{A_j(x,t,V)\,:\,(x,t)\in \Pi,\, \|V\| \le \de_0, m< j\le n\right\}\le -\Lambda_0,\\ [1mm]
\inf\left\{|A_j(x,t,V)-A_k(x,t,V)|\,:\,(x,t)\in\Pi,\, \|V\| \le \de_0, 1\le j\ne k\le n\right\}\ge \Lambda_0,\\ [1mm]
\inf\left\{\left|\det Q(x,t,V)\right|\,:\,(x,t)\in \Pi,\, \|V\| \le \de_0\right\}\ge \Lambda_0. 
\end{array}
$$
\end{itemize}

\item[\bf(A2)]
  $f\in BC^2_t(\Pi;\R^n)$,  $\d_xf\in BC^1_t(\Pi;\R^n)$,  and    $h\in BC^2(\R;\R^n)$.
\item[\bf(A3)]
$R$ is  a bounded linear operator  on $BC(\R;\R^n)$.
The restriction of $R$ to  $BC^1(\R;\R^n)$  (respectively, to $BC^2(\R;\R^n)$)
 is a bounded linear operator on  $BC^1(\R;\R^n)$
(respectively, on  $BC^2(\R;\R^n)$). Moreover, for $v\in BC^1(\R;\R^n)$ it holds
\beq\label{R'}
\begin{array}{rcl}
\displaystyle  \frac{d}{dt}(Rv)_j(t)&=& \displaystyle\left(R^\prime v\right)_j(t)
                           + \bigl(\widetilde R v^\prime\bigr)_j(t),\\ [2mm]
\displaystyle \frac{d}{dt}(\widetilde Rv)_j(t)&=& \displaystyle \bigl(\widetilde R^\prime v\bigr)_j(t)
                                 + \bigl(\widehat R v^\prime\bigr)_j(t),
                                 \end{array}
\ee
where $v^\prime(t)=\frac{d}{dt}v(t)$ and
$R^\prime,\widetilde R,\widetilde R^\prime,\widehat R: BC(\R;\R^n)\to BC(\R;\R^n)$
are some bounded linear operators.
\end{itemize}

\paragraph{Notation and further assumptions}
Set
$$
\begin{array}{cc}
a(x,t)=\diag\left(A_1(x,t,0),\dots,A_n(x,t,0)\right),\\ [2mm]
b(x,t) =Q^{-1}(x,t,0)\Bigl(B(x,t,0)Q(x,t,0)+ \partial_t Q(x,t,0) + A(x,t,0)\partial_x Q(x,t,0)\Bigr),
\end{array}
$$
and
\begin{eqnarray*}
	\displaystyle\ga_j=\inf\limits_{x,t}\frac{b_{jj}(x,t)}{|a_j(x,t)|},  \quad
	\tilde\ga_j=\inf\limits_{x,t}\left|\frac{b_{jj}(x,t)}{a_j(x,t)}\right|,\quad
	\displaystyle\be_j=\sup\limits_{x,t}\sum_{k\not=j}\left|\frac{b_{jk}(x,t)}{a_j(x,t)}\right|.
\end{eqnarray*}

In Section \ref{linear_thm} we will consider a linearized version of the 
system \reff{eq:1}; see \reff{eq:vg2}. The characteristics of this
linear system, which we need already now, are defined as follows.
 For given $j\le n$, $x \in [0,1]$, and $t \in \R$, the $j$-th characteristic 
passing through the point $(x,t)\in \Pi$, is defined
as the solution
$
\xi\in [0,1] \mapsto \om_j(\xi)=\om_j(\xi,x,t)\in \R
$
of the initial value problem
\beq\label{char}
\partial_\xi\om_j(\xi, x,t)=\frac{1}{a_j(\xi,\om_j(\xi,x,t))},\;\;
\om_j(x,x,t)=t.
\ee
Due to the assumption {\bf(A1)}, the characteristic curve $\tau=\om_j(\xi)$ reaches the
boundary of $\Pi$ in two points with distinct ordinates. Let $x_j$
denote the abscissa of that point whose ordinate is smaller.
Specifically,
$$
x_j=\left\{
 \begin{array}{rl}
 0 &\mbox{if}\ 1\le j\le m,\\
 1 &\mbox{if}\ m<j\le n.
\end{array}
\right.
$$

Write
\begin{eqnarray}
  c_j^l(\xi,x,t)=\exp \int_x^\xi
\left[\frac{b_{jj}}{a_{j}} - l\frac{\partial_t a_{j}}{a_{j}^2}\right](\eta,\om_j(\eta))\,d\eta\label{cl}\label{dl}
\end{eqnarray}
and introduce operators
 $G_1, G_2, H_1,H_2\in \LL(BC(\R,\R^n))$ by
 \beq \label{Ctilde1}
\begin{array}{rcll}
  \displaystyle \left(G_1\psi\right)_j(t) &=& c_j^1(x_j,1-x_j,t)(\widetilde R\psi)_j(\omega_j(x_j,1-x_j,t)), & \\ [1mm]
 \displaystyle \left(G_2\psi\right)_j(t) &=& c_j^2(x_j,1-x_j,t)(\widehat R\psi)_j(\omega_j(x_j,1-x_j,t)), & \\ [1mm]
\bigl(H_{l}\psi\bigr)_j(t)&=&c_j^l(x_j,1-x_j,t)\psi_j(\om_j(x_j,1-x_j,t))\quad  \mbox{if }\,\, \inf\limits_{x,t}b_{jj}>0,&\\ [1mm]
\bigl(H_{l}\psi\bigr)_j(t)&=&c_j^l(1-x_j,x_j,t)\psi_j(\om_j(1-x_j,x_j,t)) \quad  \mbox{if }\,\, \sup\limits_{x,t}b_{jj}<0.&
 \end{array}
\ee
In what follows, we will  use the simplified notations
$$
\begin{array}{rcl}
 \|R_j\| &=&\sup\left\{\|(Ru)_j\|_{BC(\R)}\,:\, \|u\|_{BC(\Pi;\R^n)}=1\right\},\\ [3mm]
 \|R\|& =&\|R\|_{\LL(BC(\R;\R^n))}=\max\limits_{j\le n}\|R_j\|.
 \end{array}
$$

We consider two sets of stable conditions on the data of the original  problem.
\begin{itemize}
	\item[\bf(B1)]
For each $j \le n$, it holds
\begin{eqnarray*} 
\begin{array}{rcl}
\|R_j\| +\frac{\be_j}{\ga_j}\left(1-e^{-\ga_j}\right)
< 1 \quad \mbox{ if }
\inf\limits_{x,t}b_{jj}> 0, \\
e^{-\ga_j}\|R_j\| +\frac{\be_j}{\ga_j}\left(1-e^{-\ga_j}\right)
< 1 \quad \mbox{ if }
\inf\limits_{x,t}b_{jj}< 0,\\
\|R_j\| +\be_j<1\quad \mbox{ if }
\inf\limits_{x,t}b_{jj}=0.
\end{array}
\end{eqnarray*}
	\item[\bf(B2)]
For each $j \le n$,	it holds
\begin{eqnarray*}
\inf\limits_{x,t}b_{jj}>0,\quad
e^{-\ga_j} \|R_j\| <1,\\
\left(1+\|R\|\left[1-\max_{i\le n} \left\{e^{-\ga_i}\|R_i\|\right\}\right]^{-1}\right)
\frac{\be_j}{\ga_j}\left(1-e^{-\ga_j}\right)<1.
\end{eqnarray*}
	\end{itemize}

Moreover, in the particular case of  periodic  boundary conditions $(Rz)_j=z_j$ or, the same, in the case
\beq\label{eq:2per}
\begin{array}{ll}
   u_{j}(0,t)=u_{j}(1,t) \,\,\mbox{ for all } j\le n,
\end{array}
\ee
we consider yet another set of conditions.
\begin{itemize}
	\item[\bf(B3)]
	For each $j \le n$, it holds 
\begin{eqnarray*}
\inf\limits_{x,t}|b_{jj}|\ne 0
\quad \mbox{and}\quad 
 \frac{\beta_j}{\tilde\ga_j}\left(2-e^{-\tilde\ga_j}\right)<1.
\end{eqnarray*}
	\end{itemize}

Note that, if $\inf\limits_{x,t}b_{jj}>0$, then
	the conditions  ${\bf (B1)}$ and  ${\bf (B2)}$
	differ at least in the restrictions imposed on the boundary operator
	$R$. More precisely, since
	the constants $\ga_j$ are positive for all $j\le n$, the condition
	${\bf (B2)}$
	allows for $\|R_j\|\ge 1$, what is not allowed
	by ${\bf (B1)}$.

\subsection{Main result}

 A continuous function $w(x,t,v)$ defined on $ [0,1]\times\R\times[-\de_0,\de_0]^n $
is  a {\it Bohr almost periodic in $t$ uniformly in $x$ and $v$ }  (see  \cite[p. 55]{Cord}) 
 if for
every $\mu > 0$ there exists a relatively dense set of $\mu$-almost
periods of~$w$, i.e., for every $\mu > 0$ there exists a positive number $l$
such that every interval of length $l$ on $\R$ contains a number $\tau$ such that
$$
| w(x,t + \tau,v) - w(x,t,v)| < \mu \ \ \mbox{ for  all }  (x,t) \in \Pi \mbox{ and } \|v\|\le\de_0.
$$
Let $AP(\R, \mathbb R^n)$ be the space of  Bohr  almost periodic vector-functions. 
 Analogously,  let $AP(\Pi, \mathbb R^n)$ 
 be the space of 
  Bohr  almost periodic vector-functions in $t$ uniformly in $x$.
By $C_T(\R, \mathbb R^n)$ and $C_T(\Pi, \mathbb R^n)$ we denote the
   spaces of continuous and $T$-periodic  in $t$  vector-functions, defined on $\R$
   and  $\Pi$, respectively.
   
  The main result of the paper is stated in the following two  theorems.
  
  \begin{thm}\label{main}
  Let the
  	conditions $\bf(A1)$--$\bf(A3)$ and at least one of  
  	the conditions ${\bf (B1)}$ and ${\bf(B2)}$ be fulfilled.
  	If the inequalities
  	\begin{eqnarray}  \label{f15}
  	\|G_i\|_{\mathcal L(BC(\R, \mathbb R^n))} < 1
  	\end{eqnarray}
  	are satisfied for both  $i=1$ and $i=2$, then the following is true:
  	
  ${\bf 1.} $
  	There exist $\eps>0$ and $\delta>0$ such that, 
  	if $\|f\|_{BC_t^2} + \|h\|_{BC^2}\le\eps$, then  
  	the problem
  	(\ref{eq:1})--\reff{z*}   has a unique classical solution $V^*\in BC^2(\Pi,\R^n)$ 
  	such that $\|V^*\|_{BC^2_t}+\|\d_xV^*\|_{BC^1_t}\le~\delta$.
  	Furthermore, there exist $\eps>0$ and $\delta>0$ such that, if  $\|f\|_{BC_t^2} +\|\d_xf\|_{BC} + \|h\|_{BC^2}\le\eps$, then  
  	$\|V^*\|_{BC^2}\le~\delta$.
  	
 ${\bf 2.} $
  	Suppose that   $A, B, Q, f$, and $h$ are Bohr almost periodic in $t$ uniformly in $x\in[0,1]$ and $V\in[-\de_0,\de_0]^n$
  	(resp., $ A, B, Q,f$, and $h$ are $T$-periodic in $t$).
  	Moreover, suppose that
  	the restriction of the boundary operator $R$ to  $AP(\R;\R^n)$  (resp., to
  	$C_T(\R, \mathbb R^n)$)
  	is a bounded linear operator on $AP(\R;\R^n)$  (resp., on
  	$C_T(\R, \mathbb R^n)$).
  	Then the bounded classical solution $V^*$ to
  	the  problem  (\ref{eq:1})--\reff{z*}  belongs to 
  	 $AP(\Pi, \mathbb R^n)$ (resp., to $C_T(\Pi, \mathbb R^n)$).
  \end{thm}
  
  \begin{thm}\label{mainH}
  	Let $(Rz)_j=z_j$ for each $j\le n$ and 
  the 	conditions $\bf(A1)$, $\bf(A2)$, and   {\bf(B3)} be fulfilled.
  	If the inequalities
  	\begin{eqnarray}  \label{f15H}
  	\|H_i\|_{\mathcal L(BC(\R, \mathbb R^n))} < 1
  	\end{eqnarray}
  	are satisfied for both $i=1$ and $i=2$,
  	then Parts 1 and 2 of Theorem \ref{main}
  	are true for the problem  (\ref{eq:1}),
  	\reff{eq:2per}, \reff{z*}.
  \end{thm}
   
The paper is organized as follows.
In Section \ref{linear_thm} we formulate statements of independent interest for general linear first order 
nonautonomous boundary value problems related to solving
the original quasilinear problem.
In Section \ref{sec:remarks} we comment on the problem 
 (\ref{eq:1})--\reff{z*} and on
 our main assumptions.
In particular, we give an example showing that in the nonautonomous 
setting the conditions \reff{f15}  and \reff{f15H} are essential
for  $C^2$-regularity of the bounded continuous solutions.
Section \ref{sec:cont} is devoted to  bounded continuous solvability of the linear boundary value problems
(including the
linearized version of the original problem).
In Section \ref{sec:regul} we prove $C^2$-regularity of the bounded continuous solutions.
A crucial point in our approach is a perturbation theorem for the general  linear
problem  \reff{eq:vg*}, \reff{eq:2}, (\ref{zam}). This result, Theorem \ref{prop3},
is proved in Section~\ref{sec:smooth_dep}.
Our main result, Theorems~\ref{main} and \ref{mainH}, is proved in 
Section~\ref{sec:proofth3a}.

\section{Relevant linear problems}\label{linear_thm}
\setcounter{equation}{0}

\paragraph{Setting}
Our approach to the quasilinear problem \reff{eq:1}--\reff{z*}
is based on a thorough analysis of a linearized problem.
As we will see later, the main reason behind global classical solvability  of the quasilinear problem
\reff{eq:1}--(\ref{z*})  lies in the fact that
the corresponding nonhomogeneous linear problem has a unique smooth bounded solution for
any smooth right-hand side. We therefore first  establish stable sufficient  conditions  ensuring  the last property.
To this end, consider the
following general nonhomogeneous linear system
\begin{equation}\label{eq:vg*}
\partial_tv  + a^*(x,t)\partial_x v + b^*(x,t) v = g(x,t), \quad x\in(0,1), t\in\R,
\end{equation}
where $g = (g_1,\dots,g_n)$ is a vector of real-valued functions,
$a^* = (a_{jk}^*)$ and $b^* =(b_{jk}^*)$ are $n\times n$-matrices of real-valued functions.
Note that, if $a^*(x,t)=A(x,t,0)$ and $b^*(x,t)=B(x,t,0)$, then \reff{eq:vg*}
is a nonhomogeneous version of the linearized system \reff{eq:1}.
This is a reason why  we use the same  notation for the general  linear problem and for the linearized version of the original quasilinear problem.

Suppose that
\beq\label{reg0}
a_{jk}^* \in BC^1(\Pi) \mbox{ and }  b_{jk}^*\in BC(\Pi)\quad \mbox{for all } j,k\le n
\ee 
and the matrix $a^*$ has $n$ real eigenvalues $a_1(x,t),\dots,a_n(x,t)$
such that 
$a_1(x,t)>\dots>a_m(x,t)> 0>a_{m+1}(x,t)>\dots>a_n(x,t)$.
Let $q(x,t) =\left (q_{jk}(x,t)\right)$  be
a nondegenerate $n\times n$-matrix 
such that $q_{jk} \in BC^1(\Pi)$ and
\begin{eqnarray}\label{eq:h111}
a(x,t)=q^{-1}(x,t) a^*(x, t)q(x,t) = \diag(a_1(x,t),\dots,a_n(x,t)).
\end{eqnarray}
The existence of such a matrix follows from the assumptions on $a^*$.
Note that, if \reff{eq:vg*} is a linearized version of \reff{eq:1}, then
the matrix $q$ is defined by $q(x,t)=Q(x,t,0)$.
Let  $\la_0$ be a positive real such that
\begin{equation}\label{eq:h11}
\begin{array}{ll}
\inf\left\{a_j(x,t)\,:\,(x,t)\in \Pi,  1\le j\le m\right\}> \lambda_0,\\ [1mm]
\sup\left\{a_j(x,t)\,:\,(x,t)\in \Pi,  m< j\le n\right\}< -\lambda_0,\\ [1mm]
\inf\left\{|a_j(x,t)-a_k(x,t)|\,:\,(x,t)\in\Pi,  1\le j\ne k\le n\right\}> \lambda_0,\\ [1mm]
\inf\left\{\left|\det q(x,t)\right|\,:\,(x,t)\in \Pi\right\}>\la_0.
\end{array}
\end{equation}

We subject the system  (\ref{eq:vg*}) to the boundary conditions
\beq\label{eq:2}
\begin{array}{cc}
	u_{j}(0,t)= (Rz)_j(t) + h_j(t), \quad 1\le j\le m,\; t\in \R, \\
	u_{j}(1,t)= (Rz)_j(t) + h_j(t), \quad m< j\le n,\; t\in \R, 
\end{array}
\ee
where
\begin{eqnarray}\label{z}
z(t)=\left(u_1(1,t),\dots,u_{m}(1,t),u_{m+1}(0,t),\dots,u_{n}(0,t)\right)
\end{eqnarray}
and
\begin{eqnarray}
\displaystyle u = q^{-1}(x,t)v. \label{zam} 
\end{eqnarray}
The system \reff{eq:vg*}  with respect to $u$ reads
\begin{equation}\label{eq:vg2}
\partial_tu  + a(x,t)\partial_x u + b(x,t) u = g(x,t), \quad x\in(0,1), \,t\in\R,
\end{equation}
where 
$b(x,t) = q^{-1}\left(b^* q + \partial_t q + a^*\partial_x q\right).$
It is evident that the  problems \reff{eq:vg*}, \reff{eq:2}, \reff{zam}
and \reff{eq:vg2},  \reff{eq:2} are equivalent.

\paragraph{An operator representation}

Let
\begin{eqnarray}
c_j(\xi,x,t)=\exp \int_x^\xi
\left[\frac{b_{jj}}{a_{j}}\right](\eta,\om_j(\eta))\,d\eta,\label{c0}\quad
d_j(\xi,x,t)=\frac{c_j(\xi,x,t)}{a_j(\xi,\om_j(\xi))}.
\end{eqnarray}
Suppose that $g$ and $h$ are sufficiently smooth and bounded functions.
As usual, a function $u\in BC^1(\Pi;\R^n)$ is called
a {\it bounded classical} solution to (\ref{eq:vg2}), (\ref{eq:2})  
 if it satisfies (\ref{eq:vg2}), (\ref{eq:2}) pointwise.
Similarly, a function $v\in BC^1(\Pi;\R^n)$  is called
a {\it bounded classical} solution to the problem  (\ref{eq:vg*}), (\ref{eq:2}), (\ref{zam}) if it satisfies  (\ref{eq:vg*}), (\ref{eq:2}),  (\ref{zam}) pointwise.
It is straightforward to show that a function $u\in BC^1(\Pi;\R^n)$
is the bounded classical solution to
(\ref{eq:vg2}), (\ref{eq:2})  if and only if $u$
satisfies the  system of integral equations
\beq\label{rep}
\begin{array}{rcl}
	&& \ u_j(x,t)= c_j(x_j,x,t)(Rz)_j(\om_j(x_j)) + c_j(x_j,x,t)h_j(\om_j(x_j))\nonumber\\ 
	  &&\quad -\displaystyle\int_{x_j}^x d_j(\xi,x,t)\left(
	\sum_{k\not=j} b_{jk}(\xi,\om_j(\xi))
	u_k(\xi,\om_j(\xi)) - g_j(\xi,\om_j(\xi))\right)d\xi,\quad  j\le n,
\end{array}
\ee
pointwise. This motivates the following definitions.
A function $u\in BC(\Pi;\R^n)$ 
is called a {\it bounded continuous} solution to (\ref{eq:vg2}), (\ref{eq:2})  if it satisfies (\ref{rep}) pointwise.
A function  $v\in BC(\Pi;\R^n)$)
is called a {\it bounded continuous} solution to (\ref{eq:vg*}), (\ref{eq:2}), (\ref{zam}) if the function $u=q^{-1}v$ satisfies (\ref{rep}) pointwise.

Let us introduce  operators
$C,D \in  \LL(BC(\Pi;\R^n))$ and an operator $F \in\LL\left(BC(\Pi;\R^{2n});BC(\Pi;\R^n)\right)$ by
\beq\label{CDF}
\begin{array}{ll}
 (Cu)_j(x,t)= c_j(x_j,x,t)(R z)_j(\om_j(x_j,x,t)),
\\ [2mm] \displaystyle
 (Du)_j(x,t)=
-\int_{x_j}^{x}d_j(\xi,x,t)\sum_{k\neq j} b_{jk}(\xi,\om_j(\xi,x,t))u_k(\xi,\om_j(\xi,x,t)) d\xi,\label{CDF}\\ \displaystyle
 \left(F(g,h)\right)_j(x,t)=\int_{x_j}^{x}d_j(\xi,x,t)g_j(\xi,\om_j(\xi,x,t)) d\xi + c_j(x_j,x,t)h_j(\om_j(x_j,x,t)).
\end{array}
\ee
Then the system \reff{rep} can be written in the operator form
\begin{eqnarray} \label{oper}
u=Cu+Du+F(g,h).
\end{eqnarray}

\paragraph{$BC$-solutions}

Theorems \ref{bounded_sol1} and \ref{bounded_sol3}
below give
{\it stable} sufficient conditions  for $BC$-solvability of the
linear problem  \reff{eq:vg*}, \reff{eq:2}, \reff{zam}. If the data 
of the problem are sufficiently smooth,  in the autonomous case 
these conditions even ensure $BC^2$-regularity.
In the nonautonomous case,  we need an additional condition to ensure $BC^1$-regularity and 
yet another condition to ensure  $BC^2$-regularity.
These additional conditions, which are stated in Theorem~ \ref{smooth_sol},
turn out to be essential; see Subsection \ref{sec:ex}.
This seems to be a new interesting phenomenon for nonautonomous hyperbolic PDEs.

\begin{thm}\label{bounded_sol1} Let $R\in\LL(BC(\R;\R^n))$.
Suppose that  the conditions  \reff{reg0}--\reff{eq:h111} and  one of  
the conditions
${\bf (B1)}$ and ${\bf(B2)}$ are fulfilled.
	Then, for any $g\in BC(\Pi;\R^n)$ and $h\in BC(\R;\R^n)$,
	the problem \reff{eq:vg*}, \reff{eq:2}, \reff{zam} has a unique bounded continuous solution $v$. Moreover, the  apriori estimate 
	\begin{eqnarray} \label{ots2}
	\|v\|_{BC}
	\le K (\|g\|_{BC} + \|h\|_{BC})
	\end{eqnarray}
is fulfilled	for a  constant $K>0$ not depending on $g$ and $h.$
\end{thm}

\begin{thm} \label{bounded_sol3}Let $(Rz)_j=z_j$ for each $j\le n$.
Suppose that  the conditions  \reff{reg0}--\reff{eq:h111} and  ${\bf(B3)}$ are fulfilled.
Then, for any $g\in BC(\Pi;\R^n)$,
the problem \reff{eq:vg*}, \reff{eq:2per}, \reff{zam} has a unique bounded continuous solution $v$. Moreover,  the estimate (\ref{ots2}) 
is fulfilled with $h=0$
and with a positive constant $K$ not depending on  $g$.
\end{thm}

Higher regularity of bounded continuous solutions is the subject of the next theorem.

\begin{thm}\label{smooth_sol}
	Assume that the assumptions of  Theorem \ref{bounded_sol1} 
	(resp., Theorem~\ref{bounded_sol3})
	are fulfilled. 
	
	${\bf 1.}$ Let $a_{jk}^*,q_{jk}\in BC^2_t(\Pi)$, $\partial_x a_{jk}^*,\partial_x q_{jk},b_{jk}^*,g_j\in BC^1_t(\Pi)$,   
	and $h_j\in BC^1(\R)$ for all $j,k\le n$. Suppose that    the restriction of $R$ to  $BC^1(\R;\R^n)$
	is a bounded linear operator on  $BC^1(\R;\R^n)$ satisfying \reff{R'}. 
	If the inequality \reff{f15} for $i=1$ (resp., the inequality  \reff{f15H} for $i=1$)
	is true, 
	then the bounded continuous solution $v$ to
	the problem \reff{eq:vg*}, \reff{eq:2},  \reff{zam}
	(resp.,  to the problem \reff{eq:vg*}, \reff{eq:2per},  \reff{zam})
	belongs to $BC^1(\Pi,\R^n)$. Moreover, the  apriori estimate 
	\begin{eqnarray} \label{ots55}
	\|v\|_{BC^1}
	\le K_1 \bigl( \| g\|_{BC_t^1} + \|h\|_{BC^1}\bigr)\ \
	\bigl(\mbox{resp., } \|v\|_{BC^1}
	\le K_1  \| g\|_{BC_t^1}\bigr)
	\end{eqnarray}
is fulfilled	for a constant $K_1>0$ not depending on $g$ and $h.$
	
	${\bf 2.}$ Let, additionally, $b^*_{jk},g_j\in BC^2_t(\Pi)$ 
	and $h_j\in BC^2(\R)$
	for all $j,k\le n$ and
	the restriction of $R$ to   $BC^2(\R;\R^n)$
	be a bounded linear operator on    $BC^2(\R;\R^n)$.
	If the inequality \reff{f15} for $i=2$ (resp., the inequality  \reff{f15H} for $i=2$)
	is true, then $v\in BC_t^2(\Pi,\mathbb R^n)$ and $\d_{x}v \in BC^1_t(\Pi,\mathbb R^n)$.
		Moreover, the  apriori estimate
	\beq \label{ots56}
		\begin{array} {cc}
	\|v\|_{BC_t^2} +\|\partial_{x} v\|_{BC_t^1}\le K_2 \bigl(\| g\|_{BC_t^2} + \|h\|_{BC^2}\bigr)\\ [2mm]
	  \bigl(\mbox{resp., } \|v\|_{BC_t^2} +\|\partial_{x} v\|_{BC_t^1}\le K_2 \| g\|_{BC_t^2}\bigr)
	\end{array}
	\ee
 is fulfilled	for a constant $K_2>0$  not depending on $g$ and $h.$
\end{thm}

\paragraph{A perturbation theorem}

Along with the system (\ref{eq:vg*}) we consider its  perturbed version
\begin{equation}\label{eq:1p}
\partial_tv  + \tilde a^*(x,t)\partial_x v + \tilde b^*(x,t) v = g(x,t), 
\end{equation}
where
$\tilde a^*=\left(\tilde a^*_{jk}\right)$ and $\tilde b^*=\left(\tilde b^*_{jk}\right)$
are $n\times n$-matrices of real-valued functions
with $\tilde a_{jk}^*\in BC^1(\Pi)$  and $\tilde b_{jk}^*\in BC(\Pi)$.
 The matrix $\tilde a^*$ is supposed to have $n$ real eigenvalues $\tilde a_1(x,t),\dots,\allowbreak\tilde a_n(x,t)$ such that 
 $\tilde a_1(x,t)>\dots>\tilde a_m(x,t)> 0>\tilde a_{m+1}(x,t)>\dots>\tilde a_n(x,t)$.
Let  $\tilde q(x,t) = \left(\tilde q_{jk}(x,t)\right)$  be a non-degenerate $n\times n$-matrix such that $\tilde q_{jk}\in BC^1(\Pi) $ and
\begin{eqnarray}\label{tildeaa}
	\tilde a(x,t)=\tilde q^{-1}(x,t)\tilde  a^*(x, t)\tilde q(x,t) = \diag(\tilde a_1(x,t),\dots,\tilde a_n(x,t)).
\end{eqnarray}
Due to the assumptions on $a^*$, we can fix $\eps_0>0$ such that, whenever $\|\tilde a^* - a^*\|_{BC}\le \eps_0$ and $\|\tilde q - q\|_{BC}\le \eps_0$,
the condition  \reff{eq:h11} is fulfilled
with $\tilde a_j$ and  $\tilde q$ in place of $a_j$ and  $q$,
respectively.

\begin{thm} \label{prop3}\
${\bf 1.}$
	If the assumptions of Part 1 of Theorem \ref{smooth_sol}  are fulfilled, then
	there exists $\eps_1 \le \eps_0$ such that, for all 
	$\tilde a_{jk}^*$, $\tilde b_{jk}^*$, and $\tilde q_{jk}$ satisfying the conditions
	\begin{eqnarray}
& 	\tilde a_{jk}^*,\tilde q_{jk}\in BC^2_t(\Pi),\ 
\partial_x \tilde a_{jk}^*, \partial_x \tilde q_{jk} \in BC^1_t(\Pi),\  \tilde b_{jk}^*\in BC^1_t(\Pi),& \nonumber	\\
&  \|\tilde a^* - a^*\|_{BC^2_t} + \|\partial_x\tilde a^* - \partial_x a^*\|_{BC^1_t} \le \varepsilon_1,\
\|\tilde b^* - b^*\|_{BC^1_t} \le \varepsilon_1,&\label{pp}\\
&  \|\tilde q- q\|_{BC^2_t} + \|\partial_x\tilde q - \partial_x q\|_{BC^1_t} \le \varepsilon_1,	\nonumber
	\end{eqnarray}
	the system (\ref{eq:1p}), \reff{eq:2}, \reff{zam} with $\tilde q$ in place of
	$q$ has a unique bounded classical solution $\tilde v\in BC^1(\Pi;\R^n)$.
	Moreover, $\tilde v$  satisfies the apriori estimate \reff{ots55} with $\tilde v$ in place of $v$ 
	for a constant $K_1$  not depending on $\tilde a^*, \tilde b^*, \tilde q, g$, and $h$.
	
	${\bf 2.}$
	If the assumptions of Part 2 of Theorem \ref{smooth_sol}  are fulfilled, then
	there exists $\eps_1 \le \eps_0$ such that, for all $\tilde a_{jk}^*$ and $\tilde b_{jk}^*$ satisfying the conditions \reff{pp}
	and the stronger conditions
	\begin{eqnarray*}
		\tilde b_{jk}^*\in BC^2_t(\Pi)\quad  \mbox{and}\quad   \|\tilde b^* - b^*\|_{BC^2_t} \le \varepsilon_1,
	\end{eqnarray*}
	the system (\ref{eq:1p}), \reff{eq:2}, \reff{zam} with $\tilde q$ in place of
	$q$ has a unique bounded classical solution $\tilde v\in BC^2_t(\Pi;\R^n)$ with $\d_x\tilde v\in BC^1_t(\Pi;\R^n).$
	Moreover, $\tilde v$  satisfies the apriori estimate \reff{ots56} with $\tilde v$ in place of $v$ for
	a constant $K_2$  not depending on $\tilde a^*, \tilde b^*, \tilde q, g$, and $h$.
\end{thm}

\section{Comments on the problem and the assumptions}\label{sec:remarks}
\setcounter{equation}{0}

\subsection{About the quasilinear system  \reff{eq:1}}
It is well known that quasilinear hyperbolic PDEs are accompanied by various singularities
 as shocks and blow-ups. Since the characteristic curves are controlled by unknown functions,
the characteristics of the same family  intersect each other  in general and, therefore,
bring different values of the corresponding unknown functions  into the  intersection points (appearance of shocks).
The nonlinearities in $B(x,t,u)$  often lead to infinite increase of solutions in a finite time (appearance of blow-ups).
When speaking about global classical solutions, one needs to provide conditions preventing the singular
behavior.

Certain classes of nonlinearities
ensuring  a non-singular behavior for autonomous quasilinear systems are described in \cite{HanNat,Li}. Some monotonicity  and sign preserving conditions on the coefficients of the nonautonomous
quasilinear hyperbolic systems are imposed
in \cite{kyr,MyFi}.
In the present paper, we study nonautonomous  quasilinear hyperbolic systems with
lower order terms and use a
different approach focusing on small solutions only. We do not need any of
the above constraints. Instead, we assume a regular behavior of the linearized system and smallness of
the right hand sides. 
Small periodic classical solutions for autonomous
quasilinear hyperbolic systems without lower order terms and with small nondiagonal elements of the matrix
$A=A(V)$ for $V\approx 0$  were investigated in  \cite{Qu}.
 In our setting, the nondiagonal entries of the matrix $A=A(x,t,V)$ are not necessarily small and the lower order coefficients $B(x,t,V)$ are 
 not necessarily zero. 
 Our dissipativity conditions depend  both on the boundary
 operator and on the coefficients of the hyperbolic system.
 
 In Section \ref{sec:ex} we show that the additional dissipativity conditions \reff{f15} and \reff{f15H}
 are in general necessary for  $C^2$-regularity of continuous solutions, which is a notable
 fact in the context of  nonautonomous hyperbolic problems.
 
\subsection{ About the boundary  conditions  \reff{eq:2*}}
The boundary operator $R$ covers different kinds of reflections and delays,
 in particular,
 $$
 (RZ)_{j}(t)=
\sum_{k=1}^{n}\left[r_{jk}(t)Z_k(t-\theta_{jk}(t))+
\int_{0}^{\vartheta_{jk}(t)}p_{jk}(t,\tau)Z_k(t-\tau)\, d\tau\right],
\quad j\le n,
$$
where $r_{jk}$,  $p_{jk}$, $\theta_{jk}$, and $\vartheta_{jk}$ are known $BC^1$-functions.  The boundary operators $R^\prime$ and $\widetilde R$ introduced  in
\reff{R'} are in this case computed by the formulas
$$
\begin{array}{rcl}
\displaystyle(R^\prime Z)_j(t)&=&\displaystyle\sum_{k=1}^{n}\Biggl[r^{\prime}_{jk}(t)Z_k(t-\theta_{jk}(t))+p_{jk}(t,\vartheta_{jk}(t))
Z_k(t-\vartheta_{jk}(t))\vartheta_{jk}^\prime(t)\\ \displaystyle &&\displaystyle+
\int_{0}^{\vartheta_{jk}(t)}\d_tp_{jk}(t,\tau)Z_k(t-\tau)\, d\tau\Biggr],\\ [2mm]
\displaystyle(\widetilde R Z)_j(t)&=&\displaystyle\sum_{k=1}^{n}\Biggl[r_{jk}(t)Z_k(t-\theta_{jk}(t))(1-\theta_{jk}^{\prime}(t))+
\int_{0}^{\vartheta_{jk}(t)}p_{jk}(t,\tau)Z_k(t-\tau)\, d\tau\Biggr].
\end{array}
$$

Boundary conditions of the reflection type appear, in particular, in semiconductor laser modeling
\cite{LiRaRe,Si} and  in boundary feedback control problems \cite{Bas,CBdA,Gug,Pav}, while integral boundary
conditions (with delays \cite{Lyu}) appear, for instance,  in hyperbolic age-structured models
\cite{car,key}.

\subsection{
	Nilpotency of the operator $C$}
Theorem \ref{bounded_sol1} can be extended
if the operator $C$ is nilpotent.
 This  is the  case of the so-called smoothing
boundary conditions, see e.g. \cite{kmit}. The smoothing property allowed us in \cite{KRT2}
to solve the  problem \reff{eq:1}--\reff{z*}
where the boundary conditions \reff{eq:2*} are specified to be of the reflection type,
without the requirement of the smallness of  $\|D\|_{\LL(BC(\Pi,\R^n))}$.
In  \cite{KRT2} we used the assumption that
  the  evolution  family  generated  by  a  linearized  problem  has  exponential
dichotomy on $\R$ and  proved that the dichotomy survives under small perturbations in the
coefficients of the hyperbolic system.
For more general boundary conditions (in particular, for \reff{eq:2*}) when the operator $C$ is
not nilpotent,  the issue of the robustness of exponential dichotomy for hyperbolic PDEs remains a challenging open problem.

\subsection{Space-periodic problems and   exponential dichotomy}
In the case of space-periodic boundary conditions \reff{eq:2per},
 our  assumption   ${\bf (B3)}$  implies, according to
\cite{KKR},  that
the  evolution  family  generated  by  the  linearized   problem
has  the exponential
dichotomy on $\R$.
For  more general boundary conditions \reff{eq:2} one can  expect the same dichotomy behavior of
the  evolution  family whenever one of the  two assumptions
 ${\bf (B1)}$ and  ${\bf (B2)}$
is fulfilled, but this still remains a subject of future work.

\subsection{Time-periodic problems and small divisors}
Analysis of time-periodic solutions to hyperbolic PDEs usually encounters a complication
known as the problem of small divisors. However, this obstacle  does not appear in our setting due to the
non-resonance assumptions  ${\bf (B1)}$,  ${\bf (B2)}$, and ${\bf (B3)}$.
Similar conditions were discussed in
\cite{kmkl,kmre}.

The completely resonant case (closely related to small divisors) is qualitatively different.
This case is discussed  in  a series of papers by Temple and Young (see, e.g., \cite{Temple1,Temple2})
about time-periodic solutions to  one-dimensional linear 
Euler equations with the periodic
boundary conditions \reff{eq:2per}.
 In this case one cannot expect any stable non-resonant conditions of
our type. More precisely, in the setting of \cite{Temple1,Temple2} it holds 
$b_{jj}=0$ for all $j$ and, hence,
our condition ${\bf (B3)}$
is not satisfied. Therefore, the operator $I-C$ is not bijective, while the bijectivity property  is a 
crucial point in our Theorems \ref{bounded_sol1} and \ref{bounded_sol3}.

\subsection{Conditions  \reff{f15} and  \reff{f15H} are essential  for higher regularity of 
bounded	continuous solutions,
in general}\label{sec:ex}
In the autonomous case, when the operator $R$ and the coefficients in the hyperbolic system
\reff{eq:vg2}  do not depend on $t$, we have $R^\prime=0$, $\widetilde R=R$, and $c_j^l\equiv c_j$ for all $j\le n$ and $l=1,2$. Then the bounds
\reff{f15}  and \reff{f15H} straightforwardly follow from the
assumptions of any of Theorems~\ref{bounded_sol1} and \ref{bounded_sol3}.
The
higher regularity of  solutions follows automatically.
This means that we have to explicitly impose the conditions \reff{f15}  and \reff{f15H} only in the nonautonomous case.

We now show that in the nonautonomous case, if the estimate
 \reff{f15} is not fulfilled for   $i=1$,  then Part 1 
of Theorem \ref{smooth_sol} is not true in general. Similarly, if \reff{f15}  is not fulfilled  for $i=2$, then one can show that Part 2 of Theorem \ref{smooth_sol} is not necessarily true.

Consider the following example,
satisfying all  the  assumptions in Part 1 of Theorem~\ref{smooth_sol} except  \reff{f15} for $i=1$:
\begin{eqnarray}
& & \partial_t u_1  + \frac{2}{4\pi-1}\partial_xu_1  = 1, \quad \partial_t u_2  - (2+\sin t)\partial_x u_2  = 0,   
\nonumber\\ 
& & u_j(x,t+2\pi) = u_j(x,t),\quad j=1,2,\label{ex}\\ 
& & u_1(0,t)=r_1(t) u_2(0,t),\quad u_2(1,t)=r_2 u_1(1,t), \nonumber
\end{eqnarray}
where $r_1$ is $2\pi$-periodic and positive $C^1$-function $r_1(t)$ and a constant $r_2$ are such that
\beq\label{rr}
 0<\sup_{t\in\R}r_1(t)<1, \quad 0 < r_2 < 1.
\ee
In this case, all assumptions of Theorem \ref{bounded_sol1} 
are true since
$\|R_1\|=\sup_{t\in\R}r_1(t)  < 1$,
$\|R_2\|=r_2<1$,
and $b_{jk}=0$ for all $j,k\le 2$.
The system \reff{ex} has a unique bounded continuous solution $u=(u_1, u_2) \in BC(\Pi,\mathbb R^2).$  Since all the 
coefficients of
the problem are $2\pi$-periodic in~$t$, it is a simple matter to show that  the solution $u$ is $2\pi$-periodic in~$t$ as well (sf.  Section~\ref{periodic}).

We have
\begin{eqnarray*}
 & & \om_1(\xi,x,t)=\frac{4\pi-1}{2}(\xi-x)+t, \\
& & \om_2(\xi,x,t)=p^{-1}(p(t)+\xi-x)  \mbox{ with } p(t)=-2t+\cos t,
\end{eqnarray*}
and
\beq\label{der}
\begin{array}{rcl}
\displaystyle\d_t\om_2(\xi,x,t)&= &\displaystyle\exp \int_\xi^x \left(\frac{a^\prime}{a^2}\right)(\om_2(\eta,x,t)) \,d \eta  \\ [3mm]
\displaystyle   &=& \displaystyle \exp \int_\xi^x\frac{d}{d\eta}\ln a(\om_2(\eta,x,t))\,d \eta\
=\frac{a(t)}{a(\om_2(\xi,x,t))} ,
\end{array}
\ee
where $a(t)=-2-\sin t$.
Then the system \reff{rep} reads
\begin{eqnarray}
& & u_1(x,t)=r_1\left( t - \frac{4\pi-1}{2}x\right)u_2\left(0, t - \frac{4\pi-1}{2}x\right)+
\frac{4\pi-1}{2}x,\label{1a} \\
& & u_2(x,t)=r_2 u_1(1,p^{-1}(p(t) + 1 - x)). \label{2a}
\end{eqnarray}

Inserting \reff{1a} into \reff{2a}, we get
\beq
\label{3a}
\begin{array}{rcl}\displaystyle
u_2(0,t)&=&\displaystyle r_2 r_1\left(p^{-1}(p(t)+1) -\frac{4\pi-1}{2}\right)\\ &&\times\displaystyle u_2\left(0, p^{-1}(p(t)+1) -\frac{4\pi-1}{2}\right) 
 + r_2 \frac{4\pi-1}{2}.
\end{array}
\ee
Using the $2\pi$-periodicity of $u_2$ in $t$, we now find the  values of $t$ at which $u_2$  has the same argument
in both sides of \reff{3a}. This is the case if, for instance,
$t - 2\pi = p^{-1}(p(t)+1) - (4\pi-1)/2$.
This equality is true if and only if $p(t)+1=p(t - \frac{1}{2})$ or, the same,
$$
\cos t-\cos \left(t - \frac{1}{2}\right)=-2\sin\left(t - \frac{1}{4}\right)\sin\left(\frac{1}{4}\right)=0.
$$
The last equation has the solutions $1/4+\pi k$, $k\in\Z$. Set $t_0=1/4$.
Then, due to \reff{rr}, the equation
 \reff{3a} yields
\beq
\label{t0}
u_2(0,t_0)=r_2\ \frac{4\pi-1}{2(1-r_2r_1(t_0))}\ne 0.
\ee
Notice the obvious identity
 $p^{-1}(p(t)+1)=\om_2(1,0,t)$. If the derivative $\d_tu_2(0,t_0)$ exists, then it is given by the formula
\beq
\label{8a}
\d_tu_2(0,t_0)=r_2 r_1(t_0)\d_t\om_2(1,0,t_0)\d_tu_2(0,t_0)+r_2r_1^\prime(t_0)\d_t\om_2(1,0,t_0)u_2(0,t_0).
\ee
By \reff{der}, we have
\begin{eqnarray*}
\d_t\om_2(1,0,t_0)=\frac{a(t_0)}{a(\om_2(1,0,t_0))}=
\frac{-2-\sin(1/4)}{-2-\sin(-1/4)}>1.
\end{eqnarray*}
 We can choose a constant $r_2$ and  a smooth $2\pi$-periodic function $r_1(t)$ such that,
additionally to the condition \reff{rr}, it holds
\beq
\label{10a}
r_2 r_1(t_0)\d_t\om_2(1,0,t_0) = 1 \mbox{ and } r_1^\prime(t_0)\not=0,
\ee
 contradicting  \reff{t0}--\reff{8a}.
This means that the continuous solution to \reff{3a} and, hence, also to \reff{1a}--\reff{2a}
 is not differentiable at $t=t_0$.
 
The violation of the condition \reff{f15}
can be seen also directly. It suffices to note that,
by~\reff{10a}, for $\psi\in BC(\R,\R^2)$ such that $\|\psi\|_{BC}=1$ and
 $\psi_1(\omega_2(1,0,t_0))\\=1$, we have
 \begin{eqnarray*}
 	\lefteqn{ \|G_1\|_{\LL(BC(\R,\R^2))}\ge |(G_1\psi)_2(t_0)|
 		= c_2^1(1,0,t_0)|(\widetilde R\psi)_2(\omega_2(1,0,t_0))|} \\ &&
 	= c_2^1(1,0,t_0)r_2|\psi_1(\omega_2(1,0,t_0))| = r_2\exp \int_0^1\left(- \frac{a^\prime(\om_2(\eta,0,t_0))}{a(\om_2(\eta,0,t_0))^2}\right)\, d\eta\\
 	& &  =  r_2\exp \int_1^0\frac{d}{d\eta}\ln a(\om_2(\eta,0,t_0))\,d\eta
 	=  r_2\frac{a(t_0)}{a(\om_2(1,0,t_0))}
 	=  r_2\d_t\om_2(1,0,t_0)>1.
 \end{eqnarray*}

\subsection{Quasilinear hyperbolic systems in applications}

Quasilinear systems of the type \reff{eq:1} cover, in particular, the one-dimensional version
of the classical Saint-Venant system for shallow water
 \cite{Sain} and  its generalizations (see, e.g. \cite{Bou}), the water flow in open-channels \cite{Hal}, and one-dimensional Euler equations  \cite{Gug2,Temple1,Temple2,Yu}.
 They are also used to describe rate-type 
materials in viscoelasticity \cite{crit,Curro,man}
and the interactions between heterogeneous cancer cell \cite{Eftimie1}.

 The behavior of  unsteady flows in  open horizontal and frictionless channel  is described in \cite{wang} by the Saint-Venant system of the type 
 \begin{eqnarray}\label{1ex}
	\partial_t A  - \partial_x (AV) = 0, \quad
	\partial_t V - \partial_x S  = 0,\qquad t\ge 0, x\in(0,L),
\end{eqnarray}
where $L$ is the length of the channel, $A=A(t,x)$ is the area of the cross section occupied by the water at position $x$ and time $t$, $V=V(t,x)$
 is the average velocity over the cross section. Furthermore,
$
S=\frac{1}{2}\,V^2+gh(A),
$
where $h(A)$ is the depth of the water. This system is 
subjected to flux boundary conditions.
Note that in the  smooth setting the system \reff{1ex}
is of our type \reff{eq:vg*}.
As   described in \cite{wang}, in a neighborhood of
an equilibrium point  the system \reff{1ex} can be written in  Riemann invariants in the diagonal  form \reff{eq:vg2}.
The flux boundary conditions are then transformed into
  boundary conditions of the type \reff{eq:2}.

The nonautonomous first order quasilinear  system 
  \begin{eqnarray*}
 \partial_t u  - \partial_x v = 0, \quad
 \partial_t v - \phi(t, v)\partial_x u  = \psi(t, v),
 \end{eqnarray*}
is  used to model the stress-strain law for metals \cite{crit,Curro,man}.
Here $v$ and $u$ denote the stress and the Lagrangian velocity, while 
the functions $\phi$ and $\psi$ measure, respectively, the non-instantaneous
 and the instantaneous response of the metal to an increment of the stress.

\section{Linear system}
\label{sec:linear}
\renewcommand{\theequation}{{\thesection}.\arabic{equation}}
\setcounter{equation}{0}

\subsection{Existence of bounded continuous solutions}\label{sec:cont}
\subsubsection{Proof of Theorem \ref{bounded_sol1}}\label{bounded_sol2}
We first give the proof under the assumption ${\bf(B1)}$.
We have to prove that  $I-C-D$ is a bijective operator from
$BC(\Pi;\R^n)$ to itself.
It suffices to establish  the estimate 
\begin{eqnarray} \label{oper0}
\|C\|_{\LL(BC(\Pi,\R^n))}+\|D\|_{\LL(BC(\Pi,\R^n))}<1.
\end{eqnarray}
Using \reff{c0},  we have
\begin{equation} \label{cj}
  \begin{array}{rcll}
\displaystyle  c_j(x_j,x,t) &=& \displaystyle\exp\left\{ \int_x^0\left[\frac{b_{jj}}{a_{j}}\right](\eta,\om_j(\eta))\,d\eta\right\}\le
e^{-\ga_j x},& j\le m,\\ [3mm]
\displaystyle \displaystyle  c_j(x_j,x,t) &=&\displaystyle \exp\left\{ \int_x^1\left[\frac{b_{jj}}{a_{j}}\right](\eta,\om_j(\eta))\,d\eta\right\}\le
e^{-\ga_j (1-x)},& j> m.
  \end{array}
  \end{equation}
  If $\inf\limits_{x,t}b_{jj}\ge 0$, then $\ga_j\ge 0$ and, if $\inf\limits_{x,t}b_{jj}< 0$, then $\ga_j< 0$.
  Combining this with
\reff{cj},  we obtain that
\begin{eqnarray*} 
	&\sup\limits_{x,t}c_j(x_j,x,t) = 1 &\mbox{if }  \ \ \inf\limits_{x,t}b_{jj}\ge 0,\\
	& \sup\limits_{x,t}c_j(x_j,x,t) \le e^{-\ga_j} &\mbox{if } \ \ \inf\limits_{x,t}b_{jj}< 0.
\end{eqnarray*}

By the definition \reff{CDF} of the operator $D$, for all $u\in BC(\Pi;\R^n)$ with $\|u\|_{BC}=1$  and 
all $(x,t)\in\Pi$ we have
\begin{eqnarray*}
  |(Du)_j(x,t)| &\le&\displaystyle
\beta_j\int_0^x\exp\left\{ \int_x^\xi
\left[\frac{b_{jj}}{a_{j}}\right](\eta,\om_j(\eta))\,d\eta\right\}\,d\xi
\le
\beta_j\int_0^xe^{-\ga_j(x-\xi)}\,d\xi\\  &=&
\frac{\be_j}{\ga_j}\left(1-e^{-\ga_j x}\right)\le \frac{\be_j}{\ga_j}
\left(1-e^{-\ga_j}\right)\qquad \mbox{if }\ j\le m,\, \ga_j\ne 0, \\[2mm]
 |(Du)_j(x,t)| &\le&\displaystyle
\beta_j\int_x^1\exp\left\{ \int_x^\xi
\left[\frac{b_{jj}}{a_{j}}\right](\eta,\om_j(\eta))\,d\eta\right\}\,d\xi \\ 
&\le& \frac{\be_j}{\ga_j}\left(1-e^{-\ga_j}\right) \qquad \mbox{if }\  j> m,\, \ga_j\ne 0, \\[2mm]
|(Du)_j(x,t)| &\le&
\be_j \qquad \mbox{if }  \ j\le n,\, \ga_j=0.
\end{eqnarray*}
Note that $\ga_j=0$ iff $\inf\limits_{x,t}b_{jj}=0$.
Using   ${\bf (B1)}$,
we immediately get   the   inequality \reff{oper0}.
This implies that, for given $g\in BC(\Pi;\R^n)$ and $h\in BC(\R;\R^n)$, the equation (\ref{oper})  has the unique solution $u=\left(I-C-D\right)^{-1}F(g,h).$
Hence, $v = q u$
is the continuous solution to (\ref{eq:vg*}), (\ref{eq:2}), (\ref{zam}). 
The estimate \reff{ots2} now easily follows.
The proof of Theorem~\ref{bounded_sol1}  under the assumption ${\bf(B1)}$ is  complete.

Now we assume that  the assumption ${\bf(B2)}$ is fulfilled. Our aim is to prove that  the operator $I-C\in\LL\left(BC(\Pi;\R^n)\right)$ is  bijective and that
the following estimate is fulfilled:
 \begin{eqnarray}\label{oper5}
	\|(I-C)^{-1}D\|_{\LL(BC(\Pi,\R^n))}<1.
\end{eqnarray}
To prove the bijectivity of  $I-C\in\LL\left(BC(\Pi;\R^n)\right)$, we  consider the equation
 \beq\label{simpl}
 u_j(x,t)=c_j(x_j,x,t)(Rz)_j(\omega_j(x_j,x,t))+r_j(x,t), \quad j\le n,
 \ee
  with respect to $u \in BC(\Pi;\R^n)$, where $r$ belongs to
  $BC(\Pi;\R^n)$ and $z$ is given by \reff{z}.
The operator  $I-C\in\LL\left(BC(\Pi;\R^n)\right)$ is bijective iff the
equation \reff{simpl} is uniquely solvable  for any $r\in BC(\Pi;\R^n)$.
Putting $x=0$ for $m<j\le n$ and $x=1$ for $1\le j\le m$,
the system \reff{simpl}
reads as follows with respect to  $z(t)$:
 \beq\label{simpl1}
   z_j(t)= c_j(x_j,1 - x_j,t)(Rz)_j(\omega_j(x_j,1 - x_j,t)) + r_j(1-x_j,t),\quad   j \le n.
 \ee
Introduce  operator $G_0 \in \LL(BC(\R,\R^n))$ by
 \beq\label{Ctilde}
 \left(G_{0} \psi\right)_j(t) = c_j(x_j,1-x_j,t)(R\psi)_j(\omega_j(x_j,1-x_j,t)), \quad   j \le n.
\ee
Note that $ \left(G_{0} z\right)_j(t) = (Cu)_j(1-x_j,t),$  $ j \le n.$
 This implies that for all $u\in BC(\Pi;\R^n)$ with $\|u\|_{BC}=1$ it holds
\begin{eqnarray} \label{--b}
 \begin{array}{ll}
\displaystyle \|(G_{0}u)_j\|_{BC} \le  \|R_j\|\exp \int_0^1
\left[-\frac{b_{jj}}{a_{j}}\right](\eta,\om_j(\eta,1,t))\,d\eta\le \|R_j\|e^{-\ga_j},& j \le m, \\[3mm]
 \displaystyle \| (G_{0}u)_j\|_{BC} \le  \|R_j\|\exp \int_0^1
\left[\frac{b_{jj}}{a_{j}}\right](\eta,\om_j(\eta,0,t))\,d\eta\le \|R_j\|e^{-\ga_j}, & j > m.
\end{array}
\end{eqnarray}
Hence, the operator $I-G_0$ is bijective  due to the assumption 
 ${\bf (B2)}$. 
 It should be noted that $\| C\|_{\LL(BC(\Pi;\R^n))}=1$, while 
$ \| G_{0}\|_{\LL(BC(\R;\R^n))}<1$.
 We, therefore, can
rewrite the system \reff{simpl1} in the form
\beq\label{z1}
z=(I-G_0)^{-1}\tilde r,
\ee
where $\tilde r(t)=\left(r_1(1,t),\dots,r_m(1,t),r_{m+1}(0,t),\dots,r_n(0,t)\right)$.
Substituting \reff{z1} into \reff{simpl}, we obtain
\beq\label{uj}
\begin{array}{rcl}
u_j(x,t)&=&\left[(I-C)^{-1}r\right]_j(x,t)\\ [2mm]
&=&c_j(x_j,x,t)\left[R(I-G_0)^{-1}\tilde r\right]_j(\omega_j(x_j,x,t))+r_j(x,t), \,\, j\le n.
\end{array}
\ee
The assumption   that $\inf\limits_{x,t}b_{jj}>0$ entails that  $c_j(x_j,x,t)\le 1$
for all $(x,t)\in\Pi$ and all $j\le n$. Therefore,
$$
\|(I-C)^{-1}\|_{\LL(BC(\Pi;\R^n))}\le \|R\|\|(I-G_0)^{-1}\|_{\LL(BC(\R;\R^n))}+1.
$$
Combining this with  the second inequality in  ${\bf (B2)}$, we arrive at the estimate
\reff{oper5}
 and, therefore,  conclude that
 the  formula 
 $u=\left[I-\left(I-C\right)^{-1}D\right]^{-1}(I - C)^{-1}F(g,h)$
gives the solution  to~(\ref{oper}).
  Hence, $v = q u$
is the continuous solution to (\ref{eq:vg*}), (\ref{eq:2}), (\ref{zam}), and this solution satisfies the  
 estimate \reff{ots2}.
 This completes the proof of Theorem \ref{bounded_sol1} 
 under the assumption ${\bf(B2)}$.
 
\subsubsection{Proof of Theorem \ref{bounded_sol3}}\label{sec:bounded_sol3}

We follow the proof of Theorem \ref{bounded_sol1}  under the assumption ${\bf(B2)}$, with the following
changes. Since in the periodic case one can integrate in both forward and backward time directions, we
use an appropriate
integral analog of the problem \reff{eq:vg2}, \reff{eq:2per}, namely
$$
\begin{array}{ll}
 u_j(x,t)= c_j(x_j,x,t)u_j(x_j,\om_j(x_j))+ c_j(x_j,x,t)h_j(\om_j(x_j))\nonumber\\
 -\displaystyle\int_{x_j}^x d_j(\xi,x,t)\left(
\sum_{k\not=j} b_{jk}(\xi,\om_j(\xi))
u_k(\xi,\om_j(\xi)) - g_j(\xi,\om_j(\xi))\right)d\xi\ \ \mbox{if } \,b_{jj}>0,\\ [6mm]
 u_j(x,t)= c_j(1-x_j,x,t)u_j(1-x_j,\om_j(1-x_j))+ c_j(1-x_j,x,t)h_j(\om_j(1-x_j))\nonumber\\
  -\displaystyle\int_{1-x_j}^x d_j(\xi,x,t)\left(
\sum_{k\not=j} b_{jk}(\xi,\om_j(\xi))
u_k(\xi,\om_j(\xi)) - g_j(\xi,\om_j(\xi))\right)d\xi\ \ \mbox{if } \, b_{jj}<0.
\end{array}
$$
Note that   in the case of general boundary conditions \reff{eq:2}
we could  integrate only in the backward time direction where the boundary
conditions are given.
Now, instead of the system \reff{simpl1}, we have  the following decoupled system:
$$
 \begin{array}{rcll}
  u_j(1-x_j,t)&=& c_j(x_j,1-x_j,t)u_j(x_j,\om_j(x_j,1-x_j,t))+r_j(1-x_j,t) \mbox{ if } \, b_{jj}>0,\\ [1mm]
  u_j(x_j,t)&= &c_j(1-x_j,x_j,t)u_j(1-x_j,\om_j(1-x_j,x_j,t))+r_j(x_j,t) \mbox{ if } \, b_{jj}<0.
 \end{array}
 $$
The analog of the operator $G_0$ introduced in \reff{Ctilde},  which we denote by $H_0$,  reads
\beq\label{Gtilde}
\begin{array}{rcll}
(H_{0}\psi)_j(t)&=&c_j(x_j,1-x_j,t)\psi_j(\om_j(x_j,1-x_j,t)) &\mbox{ if } \, b_{jj}>0,\\ [1mm]
(H_{0}\psi)_j(t)&=&c_j(1-x_j,x_j,t)\psi_j(\om_j(1-x_j,x_j,t)) &\mbox{ if } \, b_{jj}<0.
\end{array}
\ee
One can easily see that $\|C\|_{\LL(BC(\Pi;\R^n))}  = 1$, while  $\|H_{0}\|_{\LL(BC(\R;\R^n))}  < 1.$ It follows that
the operator $I-H_{0}$ and, hence, the operator $I-C$ is bijective, as desired. The rest of the proof
goes similarly to  the second part of the proof of Theorem \ref{bounded_sol1}.

\subsection{Higher regularity of the bounded continuous solutions: Proof of Theorem~\ref{smooth_sol}}\label{sec:regul}

We divide the proof into a number of claims.   Part 1 of the theorem  follows from Claims~1--4, while   Part 2   follows from Claims 5--6.

We give a proof under  the assumptions of  Theorem \ref{bounded_sol1}.
 The proof under the assumptions of  Theorem \ref{bounded_sol3}
  follows the same line, and we will
 point out only the differences. 
 
 We begin with Part 1.
Let $u \in BC(\Pi,\R^n)$ be the bounded continuous solution to the problem \reff{eq:vg2}, \reff{eq:2}.

\begin{cclaim} The generalized directional
derivatives $(\d_t + a_j\d_x)u_j$
are  continuous functions.
\end{cclaim}

\begin{subproof}
Take an arbitrary sequence of smooth functions $u^l: \Pi \to \R^n$ approaching $u$ in $BC(\Pi,\R^n)$ and an arbitrary smooth
function $\varphi: (0,1)\times\R \to \R$
with compact support. Let $\langle \cdot,\cdot\rangle$ denote the scalar product in $L^2((0,1)\times \R).$
Using \reff{rep},  for every $j\le n$ we have
\begin{eqnarray}
\lefteqn{ \left\langle(\d_t + a_j\d_x)u_j,\varphi\right\rangle = \left\langle u_j,  -\d_t \varphi - \d_x (a_j \varphi)\right\rangle =  \lim_{l \to \infty}
\left\langle u_j^l,  -\d_t \varphi - \d_x (a_j\varphi)\right\rangle}  \nonumber\\
& & =   \lim_{l \to \infty}
\Biggl\langle c_j(x_j,x,t)\left[(Rz^l)_j(\om_j(x_j,x,t)) +h_j(\om_j(x_j,x,t))\right]\nonumber\\
& & -\int_{x_j}^x d_j(\xi,x,t)\Biggl[
\sum_{k\not=j}  \left[b_{jk}u_k^l\right](\xi,\om_j(\xi)) - g_j(\xi,\om_j(\xi))\Biggr]  d\xi,
-\d_t \varphi - \d_x (a_j\varphi)\Biggr\rangle  \nonumber\\
& & = \lim_{l \to \infty}\Biggl\langle -\sum_{k=1}^n b_{jk}(x,t)u_k^l  + g_j(x,t), \varphi\Biggr\rangle
 = \Biggl\langle -\sum_{k=1}^n b_{jk}(x,t)u_k  + g_j(x,t), \varphi\Biggr\rangle, \nonumber
\end{eqnarray}
as desired. Here we used the notation 
$$
z^l(t)=\left(u_1^l(1,t),\dots,u_{m}^l(1,t),u_{m+1}^l(0,t),\dots,u_{n}^l(0,t)\right)
$$ 
and the equality
\begin{eqnarray} \label{2k}
(\d_t + a_j\d_x)\psi(\om_j(\xi,x,t)) = 0,
\end{eqnarray}
being true for  any $\psi \in C^1(\R).$
\end{subproof}

We substitute (\ref{oper}) into itself and get
\begin{eqnarray} \label{oper2}
u=Cu+(DC + D^2)u+(I + D)F(g,h).
\end{eqnarray}

\begin{cclaim} \label{2} The operators
$DC$ and $D^2$ map continuously $BC(\Pi,\R^n)$  into $BC^1_t (\Pi,\R^n)$.
\end{cclaim}
\begin{subproof}
It suffices to show that there exists a positive constant  $K_{11}$ such that for all
$u\in BC^1_t(\Pi,\R^n)$ we have
\begin{eqnarray} \label{ots3}
\left\|\partial_t\left[(DC + D^2)u\right]\right\|_{BC} \le K_{11} \|u\|_{BC}.
\end{eqnarray}
Straightforward transformations give the representation
\beq\label{ii}
\begin{array}{rcl}
  \partial_t\left[(DCu)_j(x,t)\right] &= & \displaystyle-\partial_t \Biggl(\int_{x_j}^{x}d_j(\xi,x,t)\sum_{k\neq j} b_{jk}(\xi,\om_j(\xi,x,t))
\\
& & \times c_k(x_k,\xi,\om_j(\xi))
(R z)_k(\om_k(x_k,\xi,\om_j(\xi)))\,d\xi\Biggl)  \\
&   = &- \displaystyle\sum_{k\neq j} \int_{x_j}^{x}\partial_t  d_{jk}(\xi,x,t)
(R z)_k (\om_k(x_k,\xi,\om_j(\xi)))\,d\xi  \\
& & -  \displaystyle\sum_{k\neq j} \int_{x_j}^{x}  d_{jk}(\xi,x,t)
 \frac{d}{d t}(R z)_k(\om_k(x_k,\xi,\om_j(\xi)))\,d\xi ,
\end{array}
\ee
where the functions
\begin{eqnarray*}
  d_{jk}(\xi,x,t) = d_j(\xi,x,t) b_{jk}(\xi,\om_j(\xi)) c_k(x_k,\xi,\om_j(\xi))
\end{eqnarray*}
are uniformly bounded   and have continuous and   uniformly bounded  first order derivatives in~$t$.
An upper bound as in  \reff{ots3} for the first sum in \reff{ii}
 follows directly from the regularity and the boundedness assumptions on
the coefficients of the original problem.

The strict hyperbolicity assumption \reff{eq:h11} admits the following representation formula:
\begin{eqnarray*}
 \lefteqn{\frac{d}{d t} (R z)_k(\om_k(x_k,\xi,\om_j(\xi,x,t))) } \\
 &&  = \frac{d}{d \xi} (R z)_k(\om_k(x_k,\xi,\om_j(\xi)))
 \frac{\d_3 \om_k(x_k,\xi,\om_j(\xi))\d_t\om_j(\xi)}{\d_2 \om_k(x_k,\xi,\om_j(\xi)) +
 \d_3 \om_k(x_k,\xi,\om_j(\xi))\d_\xi\om_j(\xi)}\\ [1mm]
 &&  = \frac{d}{d \xi} (R z)_k(\om_k(x_k,\xi,\om_j(\xi)))
 \frac{\d_t\om_j(\xi)a_j(\xi,\om_j(\xi))a_k(\xi,\om_j(\xi))}{a_k(\xi,\om_j(\xi))-a_j(\xi,\om_j(\xi))}.
\end{eqnarray*}
Here and in what follows $\d_j$  denotes the partial derivative with respect to the $j$-th argument.
Hence, 
\beq\label{iii}
\begin{array}{cc}
 \displaystyle\int_{x_j}^{x} d_{jk}^1(\xi,x,t) \frac{d}{d \xi} (R z)_k (\om_k(x_k,\xi,\om_j(\xi)))\,d\xi \\  \ \ \ \ =
\displaystyle d_{jk}^1(\xi,x,t) (R z)_k (\om_k(x_k,\xi,\om_j(\xi)))\Bigl|_{\xi=x_j}^{x}\\
 \displaystyle \ \ \ \ \ \ \ \ \ \ \ \ \ \ \ \  - \int_{x_j}^{x} (R z)_k (\om_k(x_k,\xi,\om_j(\xi))) \d_\xi d_{jk}^1(\xi,x,t)\, d\xi,
\end{array}
\ee
where
$$ 
d_{jk}^1(\xi,x,t) = d_{jk}(\xi,x,t)
\frac{\d_t\om_j(\xi)a_j(\xi,\om_j(\xi))a_k(\xi,\om_j(\xi))}{a_k(\xi,\om_j(\xi))-a_j(\xi,\om_j(\xi))}.
$$
Combining \reff{ii} with \reff{iii}, we conclude that  $\partial_t(DC)$ is bounded as stated in \reff{ots3}.

Similarly,
$$
\begin{array}{ll}
\displaystyle \partial_t \left[(D^2 u)_j(x,t)\right]\\
\displaystyle \quad= \sum_{k\neq j} \sum_{l\not= k} \int_{x_j}^{x} \int_{x_k}^{\xi}
\partial_t  d_{jkl}(\xi,\xi_1,x,t)
u_l(\xi_1, \om_k(\xi_1,\xi,\om_j(\xi,x,t)))\, d\xi_1 d\xi\\
\displaystyle  \quad+  \sum_{k\neq j} \sum_{l\not= k} \int_{x_j}^{x} \int_{x_k}^{\xi}d_{jkl}(\xi,\xi_1,x,t)   \d_t u_l(\xi_1, \om_k(\xi_1,\xi,\om_j(\xi,x,t)))\,d\xi_1 d\xi,
\end{array}
$$
where
$$ d_{jkl}(\xi,\xi_1,x,t) = d_j(\xi,x,t) b_{jk}(\xi,\om_j(\xi))  d_k(\xi_1, \xi, \om_j(\xi)) b_{kl}(\xi_1, \om_k(\xi_1, \xi, \om_j(\xi))).$$
A desired estimate for the first summand is obvious, while for the second summand follows from the following
transformations. For definiteness, assume that  $j ,k\le m$  (the cases $j>m$ or $k>m$ are similar).
Taking into account the identity
\begin{eqnarray*}
\lefteqn{ { \frac{d}{d t} u_l(\xi_1, \om_k(\xi_1,\xi,\om_j(\xi))) }}\\
  && = \frac{d}{d \xi} u_l(\xi_1,\om_k(\xi_1,\xi,\om_j(\xi)))
 \frac{\d_t\om_j(\xi)a_j(\xi,\om_j(\xi))a_k(\xi,\om_j(\xi))}{a_k(\xi,\om_j(\xi))-a_j(\xi,\om_j(\xi))},
\end{eqnarray*}
we get
\beq\label{ikl}
\begin{array}{ll}
 \displaystyle\int_{x_j}^{x} \int_{x_k}^{\xi} d_{jkl}(\xi,\xi_1,x,t) \d_t u_l(\xi_1, \om_k(\xi_1,\xi,\om_j(\xi)))\,d\xi_1 d\xi\\
 \displaystyle\,\,\,\,\,= \int_{x_j}^{x} \int_{x_k}^{\xi} d_{jkl}^1(\xi,\xi_1,x,t)  \frac{d}{d\xi}u_l(\xi_1, \om_k(\xi_1,\xi,\om_j(\xi)))\,d\xi_1 d\xi,
\end{array}
\ee
where
$$ d_{jkl}^1(\xi,\xi_1,x,t) =  d_{jkl}(\xi,\xi_1,x,t)
 \frac{\d_t\om_j(\xi)a_j(\xi,\om_j(\xi))a_k(\xi,\om_j(\xi))}{a_k(\xi,\om_j(\xi))-a_j(\xi,\om_j(\xi))}.
$$
The right hand side of \reff{ikl} can be written as
$$
\begin{array}{ll}
\displaystyle \int_{0}^{x} \int_{\xi_1}^{x}  d_{jkl}^1(\xi,\xi_1,x,t)
 \frac{d}{d \xi} u_l(\xi_1, \om_k(\xi_1,\xi,\om_j(\xi)))\,d\xi d\xi_1 \\ [3mm]
 \displaystyle =  \int_{0}^{x} d_{jkl}^1(\xi,\xi_1,x,t)
 u_l(\xi_1, \om_k(\xi_1,\xi,\om_j(\xi)))\Bigl|_{\xi=\xi_1}^{x}\, d\xi_1 \\ [3mm]
 \displaystyle  + \int_{0}^{x} \int_{\xi_1}^{x}
 u_l(\xi_1, \om_k(\xi_1,\xi,\om_j(\xi))) \frac{d}{d \xi} d_{jkl}^1(\xi,\xi_1,x,t) \,d\xi d\xi_1,
\end{array}
$$
which implies an upper bound as in \reff{ots3}.
\end{subproof}

\begin{cclaim}\label{3}   $I - C$ is a bijective operator from $BC^1_t (\Pi,\R^n)$ to itself.
\end{cclaim}
\begin{subproof}
We are done if we  show that the system \reff{simpl}
is uniquely solvable in $BC^1_t (\Pi,\R^n)$ for any $r\in BC^1_t (\Pi,\R^n)$.
Obviously, this  is true if and only if
\beq\label{contr1}
I-G_0 \mbox{ is a bijective operator from }   BC^1(\R,\R^n) \mbox{ to }   BC^1(\R,\R^n),
\ee
where the operator $G_0 \in \LL(BC(\R,\R^n))$ is given by \reff{Ctilde}.
To prove \reff{contr1},  let us  norm the space  $BC^1(\R,\R^n)$ with
\beq\label{beta}
\|y\|_{\sigma} = \|y\|_{BC} + \sigma \|\partial_t y\|_{BC},
\ee
where a positive constant $\sigma$ will be defined later on. Note that the norms \reff{beta} are equivalent for all $\sigma>0$.
We therefore have  to prove that  there exist  constants $\sigma<1$ and $\ga<1$ such that
$$
\|G_0 y\|_{BC}+\sigma\left\|\frac{d}{dt} G_0  y\right\|_{BC}
\le \ga\left(\|y\|_{BC} + \sigma\|y^\prime\|_{BC}\right)
\mbox{ for all } y \in  BC^1(\R,\R^n).
$$
Combining 
 \reff{c0} with the formula 
\begin{eqnarray}\label{omt}
	\partial_t \om_j(\xi,x,t)= \exp\int_\xi^x \left[\frac{\partial_ta_j}{a_j^2}\right](\eta,\om_j(\eta,x,t))\,d\eta,
\end{eqnarray}
we get that $c_j^1(\xi,x,t) = c_j(\xi,x,t)\partial_t\omega_j(\xi,x,t)$.
Then for $y\in BC^1(\R,\R^n)$ we have
\begin{eqnarray*}
& & \frac{d}{dt} (G_0 y)_j(t) =  \partial_tc_j(x_j,1 - x_j,t)(R y)_j(\om_j(x_j,1 - x_j,t)) \\
& & \hspace{25mm}   +  c_j^1(x_j,1 - x_j,t) \left[\left(R^\prime y\right)_j
+ (\widetilde R y^\prime)_j\right](\om_j(x_j,1 - x_j,t)), \  j \le n.
\end{eqnarray*}
Define an operator $W \in \mathcal L(BC(\R, \mathbb R^n))$ by
\beq\label{W}
\begin{array}{rcl}
 (W y)_j(t)& =&  \partial_tc_j(x_j,1 - x_j,t)(R y)_j(\om_j(x_j,1 - x_j,t)) \\  [2mm]&&+
c_j^1(x_j,1 - x_j,t) \left(R^\prime y\right)_j (\om_j(x_j,1 - x_j,t)).
\end{array}
\ee

On the account of \reff{Ctilde} and \reff{--b}, both  assumptions ${\bf (B1)}$   and  ${\bf (B2)}$  implies that  $\|G_0\|_{\LL(BC(\mathbb R, \mathbb R^n))}<1$. Moreover, the assumption  \reff{f15} for $i=1$
of Theorem  \ref{smooth_sol} yields  $\|G_1\|_{\LL(BC(\mathbb R, \mathbb R^n))}<1$.
Fix $\sigma<1$ such that  $\left\|G_0\right\|_{\LL(BC(\mathbb R, \mathbb R^n))}+
\sigma\left\|W\right\|_{\LL(BC(\mathbb R, \mathbb R^n))}<1$. Set
$$
\ga=\max\left\{\left\|G_0\right\|_{\LL(BC(\mathbb R, \mathbb R^n))}+
\sigma\left\|W\right\|_{\LL(BC(\mathbb R, \mathbb R^n))}, \left\|G_1\right\|_{\LL(BC(\mathbb R, \mathbb R^n))}
\right\}.
$$
Hence,
\begin{eqnarray*}
  \|G_0 y\|_{\sigma}
\le
\left\| G_0 y\right\|_{BC}
 +\sigma \| W y\|_{BC}
+ \sigma \left\|G_1 y^\prime\right\|_{BC}
\le \gamma \left(\|y\|_{BC} + \sigma\left\|y^\prime\right\|_{BC}\right).
\end{eqnarray*}
Furthermore, $\|(I - G_0)^{-1}y\|_{\sigma} \le (1 - \gamma)^{-1}\|y\|_{\sigma}$,
which yields the bound
\beq\label{ots441}
\displaystyle \|(I - G_0)^{-1}y\|_{BC_t^1} \le \displaystyle \frac{1}{\sigma}\|(I - G_0)^{-1}y\|_{\sigma} \le \displaystyle \frac{1}{\sigma(1 - \gamma)}\|y\|_{\sigma}
 \le  \frac{1}{\sigma(1 - \gamma)}\|y\|_{BC_t^1}.
\ee
Finally, from (\ref{uj}) and (\ref{ots441}) we obtain that
\beq\label{I-C--1}
\|(I-C)^{-1}\|_{\LL(BC^1_t(\Pi;\R^n))}\le 1+\frac{1}{\sigma(1 - \gamma)}\|C\|_{\LL(BC^1_t(\Pi;\R^n))}.
\ee
The proof of the claim under the assumptions of Theorem \ref{bounded_sol1} 
is complete.

The proof  under the assumptions of Theorem \ref{bounded_sol3} 
follows the same line with this changes: We specify $(Rz)_j\equiv z_j$ for all $j\le n$ and replace the operator $G_0$
by the operator $H_0$ (see the formula \reff{Gtilde}).
Hence, $(R^\prime z)_j\equiv 0$ and
$(\widetilde Rz)_j\equiv z_j$ for all $j\le n$ and all $z\in BC^1(\mathbb R, \mathbb R^n)$.
\end{subproof}

Now,   Claims \ref{2} and \ref{3} together with  the equation \reff{oper2} imply
that the bounded continuous solution $u$ to problem \reff{eq:vg2}, \reff{eq:2} is the bounded classical solution. Then, by  definition, the function 
$v=qu$
is the bounded classical solution
to the problem \reff{eq:vg*}, \reff{eq:2}, \reff{zam}.

\begin{cclaim}
	The bounded classical solution $v$ to  the problem \reff{eq:vg*}, \reff{eq:2}, \reff{zam} fulfills the  estimate \reff{ots55}.
\end{cclaim}
\begin{subproof}
Combining the estimates  (\ref{ots2}), (\ref{ots3}), and (\ref{I-C--1}) with the equations  (\ref{oper2}) and
\reff{zam}, we obtain that
$$
\begin{array}{rcl}
 \|v\|_{BC_t^1}& \le &\|q\|_{BC_t^1}\|u\|_{BC_t^1} \le\displaystyle
\|q\|_{BC_t^1}\left(1+\frac{1}{\sigma(1 - \gamma)} \|C\|_{\LL(BC^1_t(\Pi;\R^n))}\right)\\ [4mm]&\times&
  \left\|(DC + D^2)u + (I + D)F (g,h)\right\|_{BC_t^1}
  \le K_{12}\left( \| g\|_{BC_t^1}+ \| h\|_{BC^1}\right),
  \end{array}
$$
where   $u$ is the bounded classical solution to \reff{eq:vg2}, \reff{eq:2} and $K_{12}$ is a positive constant   not depending on  $g$ and $h$.
Now, from (\ref{eq:vg*})  we get
$$
\begin{array}{rcl}
 \|\partial_x v\|_{BC}&\le& \left \|(a^*)^{-1}\right\|_{BC} \left(\|g\|_{BC}
 + \|b^* v\|_{BC} + \|\partial_t v\|_{BC}\right) \\ [2mm] &\le& K_{13} \left(\| g\|_{BC_t^1}+ \| h\|_{BC^1}\right)
 \end{array}
$$
for some  $K_{13}>0$ not depending on $g$ and $h$. The estimate \reff{ots55} follows.
\end{subproof}
The proof of 
 Part 1 of the theorem  is complete.

Now we prove  Part 2.
Formal differentiation of  the system (\ref{eq:vg*}) in the distributional sense with respect to $t$ and the boundary conditions
 (\ref{eq:2})   pointwise, we get, respectively,
\begin{equation}\label{eq:112}
\begin{array}{rr}
(\partial_t + a^*\partial_x)\d_tv +
\left(b^*-\d_ta^*\,(a^*)^{-1}\right)\d_tv
+ \left(\d_tb^*-\d_ta^*\,(a^*)^{-1} \,b^*\right) v\\ [2mm]= \d_tg - \d_ta^*\,(a^*)^{-1} \,g
\end{array}
\end{equation}
and
\begin{equation} \label{eq:115}
 \d_tu_{j}(x_j,t)=\frac{d}{dt}\left(Rz\right)_j(t)+ h^\prime_j (t)=
 \left(R^\prime z\right)_j(t)+ (\widetilde Rz^\prime)_j(t) + h^\prime_j (t),
 \quad   j\le n.
\end{equation}
Introduce a new variable $w = q^{-1} \d_tv=\d_tu+q^{-1}\d_tqu$
and rewrite the problem \reff{eq:112}--\reff{eq:115} with respect to $w$
as follows:
\begin{eqnarray}\label{eq:1120}
& & \partial_t w + a(x,t)\partial_x w 
+b^1(x,t)w=g^1(x,t,v),
\end{eqnarray}
\begin{equation} \label{eq:1155}
  w_{j}(x_j,t)= \displaystyle\frac{d}{dt}\left(Rz\right)_j(t)+ h^\prime_j (t)+ \left[q^{-1}\d_tqu\right]_j(x_j,t)
  = (\widetilde Ry)_j(t) + h^1_j(t), \quad j\le n,
\end{equation}
where
 \begin{eqnarray*} 
 b^1(x,t) &=& q^{-1}\left(b^*q-\d_ta^*\,(a^*)^{-1}\,q + \d_t q + a^*\d_x q \right) =
 b - q^{-1} \d_ta^*\,(a^*)^{-1}\,q, \\ [1mm]
 g^1(x,t,v) &=& - q^{-1}\left(\d_tb^*- \d_ta^*\,(a^*)^{-1}\,b^*\right) v + q^{-1}\left(\d_tg - \d_ta^*\,(a^*)^{-1}\,g\right),\\ [1mm]
h^1_j(t) & = & \left(R^\prime z\right)_j(t)+ (\widetilde R\rho)_j(t)+ h^\prime_j (t)+ \left[q^{-1}\d_tqu\right]_j(x_j,t),
 \quad  j\le n,  \\ 
   y(t)&=&\left(w_1(1,t),\dots,w_{m}(1,t),w_{m+1}(0,t),\dots,w_{n}(0,t)\right),\\
   \rho(t)&=&\Bigl(\left[\d_tq^{-1}v\right]_1(1,t),\dots,\left[\d_tq^{-1}v\right]_m(1,t),\\               &&
           \ \   \left[\d_tq^{-1}v\right]_{m+1}(0,t),\dots,  \left[\d_tq^{-1}v\right]_{ n}(0,t)\Bigr).
 \end{eqnarray*} 
 
 \begin{cclaim}\label{cl5} The function $w\in BC(\Pi, \mathbb R^n)$ satisfies both (\ref{eq:1120}) in the distributional sense
and \reff{eq:1155} pointwise
if and only if $w$
satisfies the following system pointwise for all $j\le n$:
\beq\label{cd2}
\begin{array}{ll}
 w_j(x,t) = \displaystyle c_j^1(x_j,x,t)\left((\widetilde Ry)_j(\om_j(x_j))  + h^1_j(\om_j(x_j))\right)\\ [1mm]
\ \ \ -\displaystyle\int_{x_j}^x d_j^1(\xi,x,t)
\Biggl(\sum_{k\not=j} b_{jk}^1(\xi,\om_j(\xi)) w_k(\xi,\om_j(\xi))
 -g_{j}^1(\xi,\om_j(\xi),v(\xi,\om_j(\xi)))\Biggl)d\xi.
\end{array}
\ee
\end{cclaim}
\begin{subproof}
Set
$$
d_j^1(\xi,x,t)=\frac{c_j^1(\xi,x,t)}{a_j(\xi,\om_j(\xi))}.
$$	
 To prove the {\it sufficiency}, take an arbitrary sequence $w^l \in BC^1(\Pi;\R^n)$
 approaching $w$  in $BC(\Pi;\R^n)$. 
 Take an arbitrary smooth  function $\varphi: (0,1)\times\R \to \R$
 with compact support. On the account of \reff{cd2}, we have
\begin{eqnarray*}
& & \lefteqn{ \left\langle(\d_t + a_j\d_x)w_j,\varphi\right\rangle = -\left\langle w_j,  (\d_t \varphi + \d_x (a_j \varphi)\right\rangle  =
\lim_{l \to \infty}\left\langle w^l_j,  -\d_t \varphi - \d_x (a_j\varphi)\right\rangle }  \nonumber\\
& & =   \lim_{l \to \infty}
\Biggl\langle -c_j^1(x_j,x,t)\left((\widetilde Ry^l)_j(\om_j(x_j))  + h^1_j(\om_j(x_j))\right)\nonumber\\
& & +\int_{x_j}^x d_j^1(\xi,x,t)
\Biggl(\sum_{k\not=j} b_{jk}^1(\xi,\om_j(\xi))w_k^l(\xi,\om_j(\xi))   \\
& & - g_j^1(\xi,\om_j(\xi),v(\xi,\om_j(\xi)))\Biggr) d\xi,
\d_t \varphi + \d_x (a_j\varphi)\Biggr \rangle  \nonumber\\
& & = -\lim_{l \to \infty}\Biggl\langle \left(b_{jj}(x,t) - \frac{\partial_t a_j(x,t)}{a_j(x,t)}\right)w_j^l + \sum_{k \not= j} b_{jk}^1(x,t)w_k^l
 - g_j^1(x,t,v), \varphi\Biggr\rangle  \\
 & & = -\Biggl\langle  \left(b_{jj}(x,t) - \frac{\partial_t a_j(x,t)}{a_j(x,t)}\right)w_j + \sum_{k \not= j} b_{jk}^1(x,t)w_k
 - g_j^1(x,t,v), \varphi\Biggr\rangle,
\end{eqnarray*}
where $y^l=\left(w_1^l(1,t),\dots,w_m^l(1,t),w_{m+1}^l(0,t),\dots,w_n^l(0,t)\right)$.
It remains to note that  for all $j\le n$
 \beq\label{bjj}
  b_{jj} - a_{j}^{-1}\d_t a_j\equiv b_{jj}^1,
  \ee
what easily follows from the diagonality of the matrix $a$  and the identity
 $$
 q^{-1} \d_ta^*\,(a^*)^{-1}\,q = \left(\d_t a + q^{-1}\d_t q a - a q^{-1}\d_t q\right)a^{-1}.
 $$
Moreover, putting $x=x_j$ in \reff{cd2}, we immediately get  (\ref{eq:1155}). 
 The proof of the sufficiency is complete.
 
To prove the {\it necessity}, assume that
 the function $w$ satisfies  (\ref{eq:1120}) in the distributional sense and
(\ref{eq:1155}) pointwise.
On the account of \reff{2k},
we rewrite the system (\ref{eq:1120}) in the form
\beq\label{eq:vgc}
\begin{array}{ccc}
 (\partial_t  + a_j(x,t)\partial_x)\left(c_j^{1}(x_j,x,t)^{-1}w_j\right) \\ [3mm]=\displaystyle c_j^{1}(x_j,x,t)^{-1}
\biggl(-\sum\limits_{k\ne j}b_{jk}^1(x,t)w_k  + g_j^1(x,t,v)\biggr),
\end{array}
\ee
without destroying the equalities in the sense of distributions.
To prove that $w$ satisfies \reff{cd2}
pointwise, we use  the constancy theorem of distribution theory claiming
that any distribution on an open set with zero generalized derivatives
is a constant on any connected component of the set.
By \reff{eq:vgc}, this theorem implies that for each
$j\le n$ the expression
\beq\label{11k}
\begin{array}{ll}
 c_j^1(x_j,x,t)^{-1}\biggl[w_j(x,t)
 + \int_{x_j}^x d_j^1(\xi,x,t)\biggl(\displaystyle\sum\limits_{k\ne j}\left[b_{jk}^1w_k\right] (\xi,\om_j(\xi))\\ [4mm]\displaystyle\qquad\qquad\qquad\,\,-
g_j^1\biggl(\xi,\om_j(\xi),v(\xi,\om_j(\xi))\biggr)\biggr)\,d\xi\biggr]
\end{array}
\ee
is  constant along the characteristic curve $\om_j(\xi,x,t)$. In other words, the distributional directional derivative
$(\partial_t  + a_j(x,t)\partial_x)$ of the function \reff{11k} is equal to zero. Since \reff{11k} is a  continuous function,
$c_j^1(x_j,x_j,t)=1$,
 and the trace
$w_j(x_j,t)$ is given by  (\ref{eq:1155}), 
it follows that
 $w$ satisfies the system \reff{cd2} pointwise, as desired.
\end{subproof}

\begin{cclaim}
The bounded classical solution $v$ to  the problem \reff{eq:vg*}, \reff{eq:2}, \reff{zam} fulfills the  estimate \reff{ots56}.
\end{cclaim}

\begin{subproof}
 We rewrite the system (\ref{cd2}) in the operator form
 \begin{eqnarray} \label{oper3}
  w = C_1w + D_1 w  + F_1(g^1,h^1),
 \end{eqnarray}
where 
$C_1, D_1\in \LL(BC(\Pi;\R^n))$ and  $F_1 \in \LL(BC(\Pi;\R^{2n}), BC(\Pi;\R^n))$ are  operators defined by
\beq\label{CDF1}
\begin{array}{rcl}
 (C_1w)_j(x,t)& =&
 c_j^1(x_j,x,t)(\widetilde R y)_j(\om_j(x_j)),\\ [2mm]
 (D_1w)_j(x,t) &=&\displaystyle
 -\int_{x_j}^{x}d_j^1(\xi,x,t)\sum_{k\neq j} b^1_{jk}(\xi,\om_j(\xi))w_k(\xi,\om_j(\xi)) d\xi,  \\
\left[F_1(g^1,h^1)\right]_j(x,t) &=&\displaystyle\int_{x_j}^{x}d_j^1(\xi,x,t) 
g_j^1\left(\xi,\om_j(\xi),v(\xi,\om_j(\xi))\right)\,d\xi\\ [4mm]&&+ c_j^1(x_j,x,t)h^1_j(\om_j(x_j)).
\end{array}
\ee
 Iterating (\ref{oper3}), we obtain
  \begin{eqnarray} \label{oper4}
  w = C_1w + (D_1 C_1 + D_1^2) w + (I+D_1) F_1(g^1,h^1).
 \end{eqnarray}
 Using the same argument as in Claim 2, we conclude that the
operators $D_1 C_1$ and $D_1^2$ map continuously $BC( \Pi, \mathbb R^n)$ into $BC^1_t(\Pi, \mathbb R^n).$
Moreover, the following smoothing estimate is true:
\beq\label{smooth}
\left\|(D_1C_1 + D_1^2)w\right\|_{BC^1_t} \le K_{21} \|w\|_{BC}
\ee
for some $K_{21}>0$ not depending on $w\in BC(\Pi,\R^n)$.

Next, we prove that $I - C_1$ is a bijective operator   from $BC^1_t( \Pi, \mathbb R^n)$ to itself. 
The proof is similar to the proof of Claim \ref{3}. We  have to show that the system
$$
 w_j(x,t)=c_j^1(x_j,x,t)(\widetilde Ry)_j(\omega_j(x_j,x,t))+\al_j(x,t), \quad j\le n,
 $$
is uniquely solvable in $BC^1_t (\Pi,\R^n)$ for each $\al\in BC^1_t (\Pi,\R^n)$.
It is sufficient to show that
\beq\label{contr11}
  I-G_1 \mbox{ is a bijective operator from }   BC^1(\R,\R^n) \mbox{ to itself},
\ee
where the operator $G_1 \in \LL(BC(\R,\R^n))$ is given by \reff{Ctilde1}.
To prove \reff{contr11},  we  use 
the space  $BC^1(\R,\R^n)$ normed by \reff{beta}. 
We are  done if we prove  that  there exist  constants $\sigma_1<1$ and $\ga_1<1$ such that
$$
\|G_1 \psi\|_{BC}+\sigma_1\left\|\frac{d}{dt}G_1  \psi\right\|_{BC}
\le \ga_1\left(\|\psi\|_{BC}+\sigma_1\|\psi^\prime\|_{BC}\right)
$$
for all $\psi \in  BC^1(\R,\R^n)$. 

As it follows from \reff{dl} and \reff{omt},
$c_j^2(\xi,x,t) = c_j^1(\xi,x,t)\partial_t\omega_j(\xi,x,t)$.
Define  operator $W_1 \in \mathcal L(BC(\mathbb R, \mathbb R^n))$ by
\beq\label{W1}
\begin{array}{rcl}
	  (W_1 \psi)_j(t) &=&  \partial_tc_j^1(x_j,1-x_j,t)(\widetilde R \psi)_j(\om_j(x_j,1-x_j,t))  \\ [2mm]
	& & + c_j^2(x_j,1-x_j,t) \bigl(\widetilde R^\prime \psi\bigr)_j (\om_j(x_j,1-x_j,t)), \  j \le n.
\end{array}
\ee
Taking into  account  \reff{R'} and \reff{Ctilde1}, for given $\psi\in BC^1(\R,\R^n)$,
it holds
\begin{eqnarray*}
 \frac{d}{dt} \left[(G_1\psi)_j(t)\right]& =& \partial_tc_j^1(x_j,1-x_j,t)(\widetilde R \psi)_j(\om_j(x_j,1-x_j,t)) \\ [2mm]
& &  +  c_j^2(x_j,1-x_j,t) \Bigl[(\widetilde R^\prime \psi)_j
 + (\widehat R \psi^\prime)_j\Bigr](\om_j(x_j,1-x_j,t))\\ [2mm]
 &  =&  (W_1 \psi)_j(t)+(G_2 \psi)_j(t), \quad  j \le n.
\end{eqnarray*}
By  the assumption (\ref{f15}),
  $\|G_1\|_{\LL(BC(\mathbb R, \mathbb R^n))}<1$ and  $\|G_2\|_{\LL(BC(\mathbb R, \mathbb R^n))}<1$.
Fix $\sigma_1<1$ such that  $\left\|G_1\right\|_{\LL(BC(\mathbb R, \mathbb R^n))}+
\sigma_1\left\|W_1\right\|_{\LL(BC(\mathbb R, \mathbb R^n))}<1$. Set
$$
\ga_1=\max\left\{\left\|G_1\right\|_{\LL(BC(\mathbb R, \mathbb R^n))}+
\sigma_1\left\|W_1\right\|_{\LL(BC(\mathbb R, \mathbb R^n))}, \left\|G_2\right\|_{\LL(BC(\mathbb R, \mathbb R^n))}
\right\}.
$$
It follows that
\begin{eqnarray*}
\lefteqn{ \|G_1 \psi\|_{\sigma_1} = \|G_1\psi\|_{BC} +
\sigma_1 \left\|\frac{d}{dt} G_1\psi\right\|_{BC} \le
\left\| G_1 \psi\right\|_{BC} } \\
& & + \sigma_1 \| W_1 \psi\|_{BC}
+ \sigma_1 \left\|G_2 \psi^\prime\right\|_{BC}
\le \gamma_1 \left(\|\psi\|_{BC} + \sigma_1\left\|\psi^\prime\right\|_{BC}\right),
\end{eqnarray*}
which proves \reff{contr11}.

Similarly to \reff{I-C--1}, the inverse to $I-C_1$ can be estimated from above as follows:
$$
\|(I-C_1)^{-1}\|_{\LL(BC^1_t(\Pi;\R^n))}\le 1+\frac{1}{\sigma_1(1 - \gamma_1)}\|C_1\|_{\LL(BC^1_t(\Pi;\R^n))}.
$$
Combining this estimate with \reff{ots55}, (\ref{oper4}),   and  \reff{smooth}, we get
$$
\begin{array}{ll}
  \|\partial_t v\|_{BC^1_t} \le\|q\|_{BC^1_t}\|w\|_{BC^1_t}\\ [2mm]
\le\displaystyle\|q\|_{BC^1_t} \left(1+\frac{\|C_1\|_{\LL(BC^1_t(\Pi;\R^n))}}{\sigma_1(1 - \gamma_1)}\right)
 \left\|(D_1C_1 + D_1^2)w + (I + D_1) F_1(g^1,h^1)\right\|_{BC^1_t} \\ [4mm]
                            \le K_{22} 
                             \left(\| g^1\|_{BC^1_t}
                             +\| h^1\|_{BC^1}\right)\le 
                             K_{23} \left(\| g\|_{BC_t^2}+\| h\|_{BC^2}\right),
  \end{array}
$$
the constants $K_{22} $ and $K_{23} $  being independent of $g$ and $h$.
By \reff{eq:112},  there exists a constant $K_{24} $   not depending on $g$ and $h$ such that
$$
 \|\partial_x v\|_{BC^1_t}   \le K_{24} \left(\| g\|_{BC_t^2}+\| h\|_{BC^2}\right),
$$
which  implies the estimate \reff{ots56}, as desired.
\end{subproof}

\subsection{A perturbation result: Proof of Theorem \ref{prop3}}\label{sec:smooth_dep}

In the new variable
\beq\label{zam1}
u=\tilde q^{-1}v,
\ee
the perturbed system \reff{eq:1p}
reads
\begin{equation}\label{eq:vg21}
\partial_tu  + \tilde a(x,t)\partial_x u + \tilde b(x,t) u = g(x,t), 
\end{equation}
where $\tilde a$ is defined by \reff{tildeaa} and
$\tilde b(x,t) = \tilde q^{-1}\left(\tilde b^* \tilde q + \partial_t \tilde q + \tilde a^*\partial_x \tilde q\right).$

We  will use the following notation.
 The $j$-th characteristic of \reff{eq:vg21}
passing through the point $(x,t)\in \Pi$ is defined
as the solution $\xi\in [0,1] \mapsto \tilde\om_j(\xi)=\tilde\om_j(\xi,x,t)\in \R$ of the initial value problem
$$
\partial_\xi\tilde\om_j(\xi, x,t)=\frac{1}{\tilde a_j(\xi,\tilde\om_j(\xi, x,t))},\;\;
\tilde\om_j(x,x,t)=t.
$$
Introduce notation
$$
\begin{array} {rl}
\displaystyle \tilde c_j(\xi,x,t)=\exp \int_x^\xi
\Biggl[\frac{\tilde b_{jj}}{\tilde a_{j}}\Biggr](\eta,\tilde\om_j(\eta))\,d\eta,&\displaystyle
\quad\tilde d_j(\xi,x,t)=\frac{\tilde c_j(\xi,x,t)}{\tilde a_j(\xi,\tilde\om_j(\xi))} 
\end{array}
$$
and operators
$\widetilde C, \widetilde D \in \LL(BC(\Pi;\R^n))$ and  $\widetilde F\in \LL\left(BC(\Pi;\R^{2n});BC(\Pi;\R^n)\right)$ by
\begin{eqnarray*}
(\widetilde C u)_j(x,t)&=& \tilde c_j( x_j,x,t)(R u)_j(\tilde\om_j(x_j)), \\
(\widetilde Du)_j(x,t)&=&
 -\int_{ x_j}^{x} \tilde d_j(\xi,x,t)\sum_{k\neq j} \tilde b_{jk}(\xi,\tilde\om_j(\xi))u_k(\xi,\tilde\om_j(\xi)) d\xi,\\
(\widetilde F(g,h))_j(x,t)&=&\int_{ x_j}^{x}\tilde d_j(\xi,x,t)g_j(\xi,\tilde \om_j(\xi)) d\xi+ 
\tilde c_j( x_j,x,t)h_j(\tilde\om_j(x_j)).
\end{eqnarray*}
Consider the following operator representation of  the perturbed problem \reff{eq:vg21}, \reff{eq:2}, \reff{zam1}:
\begin{eqnarray} \label{oper-p}
u=\widetilde Cu + \widetilde D u + \widetilde F(g,h).
\end{eqnarray}
Iterating  this formula, we get
\begin{eqnarray} \label{oper-pp}
u=\widetilde Cu + \left(\widetilde D\widetilde C + \widetilde D^2\right)u+\left(I+\widetilde D\right)\widetilde F(g,h).
\end{eqnarray}

Let $\widetilde G_0, \widetilde G_1, \widetilde G_2,\widetilde W\in \LL(BC(\R,\R^n))$,   $\widetilde C_1, \widetilde D_1\in \LL(BC(\Pi,\R^n))$, and   $\widetilde F_1\in \LL(BC(\Pi,\R^{2n}))$ denote operators given by the right hand sides of the formulas \reff{Ctilde},
\reff{Ctilde1},   \reff{W}, and \reff{CDF1}, respectively, where $ a$, $ b$, and
$\om_j$ are replaced by $\tilde a$, $\tilde b$,
and  $\tilde\om_j$, respectively.

Assume that the condition ${\bf (B1)}$ is fulfilled.  Similar argument works in the case of
${\bf (B2)}$ or ${\bf (B3)}$.

 {\it Proof of Part 1.}
  Note that the  assumptions   ${\bf (B1)}$ and \reff{f15} of Theorem
 \ref{smooth_sol} are stable with respect to small perturbations of $a$ and $b$. Since small perturbations of
 $a^*$, $b^*$, and $q$ imply small perturbations of $a$ and $b$, 
 there exists $\eps_{11}\le\eps_0$ such that, for all $\tilde a^*$ and $\tilde b^*$ varying in the range
 \beq\label{ppp}
 \begin{array}{cc}
 \|\tilde a^* - a^*\|_{BC^2_t} +
 \|\partial_x\tilde a^* - \partial_x a^*\|_{BC^1_t} \le \varepsilon_{11},
\quad \|\tilde b^* - b^*\|_{BC^1_t} \le \varepsilon_{11},\\ [2mm]
 \|\tilde q - q\|_{BC^2_t} +
\|\partial_x\tilde q - \partial_x q\|_{BC^1_t} \le \varepsilon_{11},
 \end{array}
 \ee
 the conditions  ${\bf (B1)}$ and \reff{f15} for $i=1$ remain to be true
 with $\tilde a$ and $\tilde b$
 in place of $a$ and $b$, respectively. 
 
 Due to Part 1 of Theorem  \ref{smooth_sol}, 
   the system (\ref{eq:vg21}), \reff{eq:2}, \reff{zam1} has a unique
bounded classical solution $\tilde u\in BC^1(\Pi;\R^n)$ for each fixed 
$\tilde a^*$, $\tilde b^*$, and $\tilde q$.

To derive the apriori estimate \reff{ots55} with $\tilde v$ in place of $v$,
note that  the value of 
 $\eps_{11}>0$  can be chosen so  small  that there exists 
 a   positive real $\nu_1<1$ such that  the
left hand sides of  ${\bf (B1)}$ and \reff{f15} for $i=1$, calculated for the perturbed problem
(\ref{eq:vg21}), \reff{zam1}, \reff{eq:2}, are bounded from above by $1-\nu_1$.
Due to the proof of Theorem~\ref{bounded_sol1}, this implies the inequality
$$
\|\widetilde C\|_{\LL(BC(\Pi;\R^n))}+\|\widetilde D\|_{\LL(BC(\Pi;\R^n))}\le 1-\nu_1,
$$
which is  uniform in $\tilde a^*$, $\tilde b^*$, and $\tilde q$.
Combining this inequality  with \reff{oper-p},
we conclude that there exists a  constant $\widetilde K>0$ not depending on $\tilde a^*, \tilde b^*$,  $\tilde q$, $g$, and $h$ such that
\begin{eqnarray} \label{estBC}
  \|\tilde u\|_{BC}  \le \widetilde  K\left(\| g\|_{BC}+\| h\|_{BC}\right).
 \end{eqnarray}
We immediately see from  ${\bf (B1)}$, \reff{f15} for $i=1$, \reff{--b}, and \reff{W} that there exist constants $\widetilde K_{1}>0$ and $\nu_2<1$
such that
\begin{eqnarray*}
	&&\|\widetilde G_0\|_{\LL(BC(\R,\R^n))}\le 1-\nu_2,\ \|\widetilde G_1\|_{\LL(BC(\R,\R^n))}\le 
	1-\nu_2,\ \|\widetilde W\|_{\LL(BC(\R,\R^n))}\le \widetilde K_{1},
\end{eqnarray*}
uniformly in
$\tilde a^*, \tilde b^*$, and $\tilde q$ fulfilling 
\reff{ppp} with  $\eps_{12}$ in place of~$\eps_{11}$.

Put $
\ga=1-\nu_2+\sigma \widetilde K_{1}
$
and  fix
 $\sigma<1$ such that  $\ga<1$.
Now we apply the  argument used to prove the estimate  \reff{I-C--1} and  get
\beq\label{I-C--11}
\|(I-\widetilde C)^{-1}\|_{\LL(BC^1_t(\Pi;\R^n))}\le 1+\frac{1}{\sigma(1 - \gamma)}\|\widetilde C\|_{\LL(BC^1_t(\Pi;\R^n))}.
\ee
Similarly to the proof of Claim 2 in Section \ref{sec:regul}, we show
 that the operators $\widetilde D\widetilde C$
 and $\widetilde D^2$ are smoothing and map
 $BC(\Pi,\R^n)$ into $BC^1_t(\Pi,\R^n)$. Moreover,   there exists a 
 constant $\widetilde K_{2}$ such that, for all 
$\tilde a^*$, $\tilde b^*$, and $\tilde q$
fulfilling the inequalities \reff{ppp} with $\varepsilon_{12}$
in place of $\varepsilon_{11}$,
it holds
\begin{eqnarray} \label{ots31}
\left\|\partial_t\left[\bigl(\widetilde D\widetilde C + \widetilde D^2\bigr)u\right]\right\|_{BC} \le \widetilde K_{2} \|u\|_{BC}
\end{eqnarray}
for all $u\in BC(\Pi,\R^n)$.
Now
we combine the estimates \reff{estBC}--\reff{ots31} with the equations \reff{oper-pp}.  We conclude that there exist constants
 $\eps_1\le \eps_{12}$
and  $K_1>0$ such that, for all  $\tilde a^*$, $\tilde b^*$, $\tilde q$, $g$, and $h$ varying in the range \reff{ppp} with $\varepsilon_{1}$
in place of $\varepsilon_{11}$, the estimate  \reff{ots55} is true with
$\tilde v$ in place of $v$.

{\it Proof of Part 2.}
Let $\eps_1$ be a constant satisfying Part 1 of Theorem \ref{prop3}.
Consider a perturbed version of the equation \reff{oper4}
where $C_1$, $D_1$, and $F_1$ are replaced by 
$\widetilde C_1$, $\widetilde D_1$, and $\widetilde F_1$, respectively.
Proceeding similarly to Part 1, we use \reff{W1} and \reff{f15} for $i=2$
and conclude that the
constant $\eps_1$ can be chosen so small that there exist  positive reals 
$\nu_3<1$ and $\widetilde K_{3}$ fulfilling the bounds
\begin{eqnarray}
\|\widetilde G_2\|_{\LL(BC(\R,\R^n))}<1-\nu_3,\quad\|\widetilde W_1\|_{\LL(BC(\R,\R^n))}\le \widetilde K_{3},\label{G2}
\end{eqnarray}
uniformly in  $\tilde a^*$, $\tilde b^*$, and $\tilde q$
satisfying the estimates \reff{ppp} with $\varepsilon_{1}$
in place of $\varepsilon_{11}$ as well as the stronger estimate
 $ \|\tilde b^* - b^*\|_{BC^2_t} \le \varepsilon_1$.
The desired  apriori estimate \reff{ots56} for the $\eps_1$-perturbed
problem then easily follows from the perturbed versions of
 \reff{oper4} and \reff{smooth}. 

  The proof of  Theorem \ref{prop3} is complete.
  
\section{Quasilinear system: Proof of main result}\label{sec:proofth3a}
\renewcommand{\theequation}{{\thesection}.\arabic{equation}}
\setcounter{equation}{0}

\subsection{Proof of Part 1 of Theorem \ref{main}: Bounded solutions}\label{sec:bounded}
Let $\de_0$ be a constant satisfying Assumption {\bf(A1)} and  $\eps_1$ be a constant satisfying
Part 2 of Theorem \ref{prop3}.
Since the functions $A$ and $B$ are $C^2$-smooth,
there exists  $\delta_1\le\de_0$ such that for all $\varphi\in BC^2( \Pi, \mathbb R^n)$
with
\beq\label{phiest}
\|\varphi\|_{BC^2_t} +\|\partial_x \varphi\|_{BC^1_t}\le \delta_1
\ee
we have
\beq\label{minus}
\| a^\varphi \|_{BC^2_t}  +  \| \partial_xa^\varphi \|_{BC^1_t} \le \eps_1, \quad \| q^\varphi \|_{BC^2_t}  +  \| \partial_xq^\varphi \|_{BC^1_t} \le \eps_1, \quad \| b^\varphi \|_{BC^2_t} \le \eps_1,
\ee
where $a^\varphi (x,t) = A(x,t,\varphi (x,t)) - A(x,t,0)$,
$q^\varphi (x,t) = Q(x,t,\varphi (x,t)) - Q(x,t,0)$,
 and $b^\varphi (x,t) = B(x,t,\varphi (x,t)) - B(x,t,0). $
Therefore,  due to Theorem \ref{prop3}, for given $\vphi\in BC^2(\Pi;\R^n)$ satisfying  \reff{phiest},
 the system
\beq\label{eq:phi}
\partial_t V + A(x,t,\varphi)\partial_x V + B(x,t, \varphi)V = f(x,t)
\ee
with the boundary conditions (\ref{eq:2*}) and with $U(x,t) =  Q^{-1}(x,t,\vphi)V(x,t)$ has a unique solution 
$V^\varphi \in BC_t^2(\Pi,\mathbb R^n)$ such that $\d_{x}V^\varphi \in BC^1_t(\Pi,\mathbb R^n)$.
Moreover,   the estimate  (\ref{ots56}) holds with $v$ replaced by
$V^\varphi$  and is uniform in $\vphi$
obeying  \reff{phiest}. Since  $V^\varphi$ satisfies \reff{eq:phi}, it belongs to $BC^2(\Pi,\mathbb R^n)$.

 Put $V^0(x,t) =0.$  For a given nonnegative integer number $k$, construct the iteration  $V^{k+1}(x,t)$ as the unique bounded classical
solution   to the linear system
\begin{eqnarray} \label{th-2}
 \partial_t V^{k+1}+ A(x,t,V^k)\partial_x V^{k+1}  + B(x,t,V^k)V^{k+1} = f(x,t)
\end{eqnarray}
subjected to the boundary conditions 
\beq\label{eq:2*k}
\begin{array}{ll}
	U^{k+1}_{j}(x_j,t)= (RZ^{k+1})_j(t) + h_j(t), \;\;\;  j\le n, 
\end{array}
\ee
where
$$
Z^{k+1}(t)=\left(U_1^{k+1}(1,t),\dots,U_{m}^{k+1}(1,t),U_{m+1}^{k+1}(0,t),\dots,U_{n}^{k+1}(0,t)\right)
$$
and
\beq\label{zamk}
U^{k+1}(x,t) =  Q^{-1}(x,t,V^k)V^{k+1}(x,t).
\ee
The function $U^{k+1}$ then satisfies the system
\begin{eqnarray*} 
 \partial_t U^{k+1}+ \hat a(x,t,V^k)\partial_x U^{k+1}  + \hat b(x,t,V^k)U^{k+1} = Q^{-1}(x,t,V^k)f(x,t),
\end{eqnarray*}
where
 \beq\label{hat}
 \begin{array}{rcl}
\hat a(x,t,V^k)&=&\diag\left(A_1(x,t,V^k),\dots,A_n(x,t,V^k)\right), \\ [1mm]
 \hat b(x,t,V^k) &=& (Q^k)^{-1}\left(B^k Q^k + \d_t Q^k + A^k\d_x Q^k \right).
\end{array}
\ee
Here and below in this proof we also use  the short notation  $A^k$, $B^k$,  $Q^k$, $a^k$, and $b^k$ for  $A(x,t,V^k)$,
$B(x,t,V^k)$,   $Q(x,t,V^k)$, $\hat a(x,t,V^k)$, and $\hat b(x,t,V^k)$, respectively.  

We divide the proof into a number of claims.

\begin{cclaim}\label{ots5}
Suppose that
\begin{eqnarray} \label{sm7aa}
\|f\|_{BC^2_t}+\|h\|_{BC^2} \le \delta_1/ K_2,
\end{eqnarray}
where $K_2$ is the constant  as in Part 2 of
Theorem \ref{prop3}.
Then there exists  a sequence  $V^k$ of bounded classical solutions to \reff{th-2}--(\ref{zamk})  belonging to $BC^2(\Pi;\R^n)$ such that
\beq\label{apr_seq}
 \| V^k\|_{BC^2_t}+ \|\partial_{x} V^k\|_{BC^1_t}\le \delta_1 \quad \mbox{for all } k.
\ee
\end{cclaim}

\begin{subproof}
  Note that the first iteration $V^1$  satisfies the system (\ref{eq:vg*}) with $g=f$ and the boundary conditions
  \reff{eq:2}, \reff{zam}.
Due to Theorem \ref{smooth_sol}, there exists a unique bounded
classical solution $V^1$ such that $V^1 \in BC_t^2(\Pi,\mathbb R^n)$ and $\d_{x}V^1 \in BC_t^1(\Pi,\mathbb R^n)$.
Since $A^0$,  $Q^0$, and $B^0$ are continuously differentiable in $x$, from  the
system  \reff{th-2} differentiated in $x$ it follows that
 $V^1 \in BC^2(\Pi,\mathbb R^n)$.
Moreover,  $V^1$  satisfies the bound (\ref{ots56}) with $v$ and $g$ replaced by  $V^1$ and $f$, respectively.
 Since  $ f$  and $h$ obey  \reff{sm7aa}, the estimate \reff{apr_seq} with $k=1$
 follows.
 Due to \reff{phiest}--\reff{minus}, we then have
 \beq\label{difference}
 \begin{array}{cc}
 \|A^1 - A^0\|_{BC^2_t} +
 \|\d_{x} A^1 - \d_{x}A^0\|_{BC^1_t}\le \varepsilon_1, \quad
 \|B^1 - B^0\|_{BC^2_t} \le \varepsilon_1,\\ [2mm]
 \|Q^1 - Q^0\|_{BC^2_t} +
 \|\d_{x} Q^1 - \d_{x}Q^0\|_{BC^1_t}\le \varepsilon_1.
\end{array}
 \ee
Theorem \ref{prop3} now implies that there exists a unique   bounded
classical solution $V^2$
such that $V^2 \in BC_t^2(\Pi,\mathbb R^n)$ and $\d_{x}V^2 \in BC^1_t(\Pi,\mathbb R^n)$.
Similarly,   $V^2(x,t)$
fulfills the bound \reff{apr_seq} with $k=2$ and, due to \reff{th-2}, belongs to $BC^2(\Pi,\mathbb R^n)$. On the account of \reff{phiest}--\reff{minus},
we also have the estimates \reff{difference} with $A^1$,
 $B^1$, and  $Q^1$ replaced by $A^2$,
 $B^2$, and  $Q^2$, respectively.

Proceeding by induction, assume that the problem \reff{th-2}--(\ref{zamk})
has  a unique   bounded
classical solution $V^k$  belonging to
$BC^2(\Pi,\mathbb R^n)$ and satisfying the estimate \reff{apr_seq} and, hence the estimates 
  \beq\label{differencek}
 \begin{array}{cc}
 	 \|A^k - A^0\|_{BC^2_t} +
 \|\d_{x} A^k - \d_{x}A^0\|_{BC^1_t}\le \varepsilon_1, \quad
 \|B^k - B^0\|_{BC^2_t} \le \varepsilon_1,\\ [2mm]
 \|Q^k - Q^0\|_{BC^2_t} +
 \|\d_{x} Q^k - \d_{x}Q^0\|_{BC^1_t}\le \varepsilon_1.
\end{array}
\ee
 Now, using  Theorem \ref{prop3} and  the system \reff{th-2}
 differentiated in $x$, we conclude that the problem \reff{th-2}--(\ref{zamk}) has  a unique   bounded
classical
 solution $V^{k+1}\in BC^2(\Pi,\mathbb R^n)$. Moreover, this solution 
   fulfills the inequalities \reff{apr_seq}
 and \reff{differencek}
with $k+1$ in place of $k$.
\end{subproof}

\begin{cclaim}\label{convV}
There exists $\varepsilon_2 \le \delta_1/ K_2 $ such that, if
$\|f\|_{BC^2_t}+\|h\|_{BC^2} < \varepsilon_2,$ then 
the sequence  $V^k$ of solutions to the problem \reff{th-2}--(\ref{zamk})
 converges in   $BC^1(\Pi;\R^n)$ to a classical solution to~(\ref{eq:1})--(\ref{z*}).
\end{cclaim}

\begin{subproof} Set 
\begin{eqnarray}
w^{k+1}& =& V^{k+1} - V^k= Q^kU^{k+1}-Q^{k-1}U^{k},\label{zam2*} \\
Y^{k+1} &=& \left(Q^k\right)^{-1}w^{k+1}.\label{zam2}
\end{eqnarray}
Hence,
 \beq\label{zam3}
 Y^{k+1} = U^{k+1} -  U^{k}  + \chi^k,
 \ee
 where $\chi^k=(Q^k)^{-1}\left(Q^k - Q^{k-1}\right) U^k$.

First we derive a boundary value problem for $w^{k+1}$. To this end, introduce the 
following notation:
$$
\begin{array}{rcl}
\bar \chi^k(t)&=&\left(\chi_1^k(1,t),\dots,\chi_{m}^k(1,t),\chi_{m+1}^k(0,t),\dots,\chi_{n}^k(0,t)\right),\\ [2mm]
\bar Y^{k+1}(t)&=&Z^{k+1}-Z^{k}+\bar \chi^k\\ &=&\Bigl(Y_1^{k+1}(1,t),\dots,Y_{m}^{k+1}(1,t),Y_{m+1}^{k+1}(0,t),\dots,Y_{n}^{k+1}(0,t)\Bigr),
\\ [2mm]
\zeta^k_j(t)&=&- \left[R\left(\bar \chi^{k}\right)\right]_j(t) 
+ \chi_j^k(x_j,t), \  j\le n.
\end{array}
$$
On the  account of \reff{zam2} and  \reff{zam3}, the boundary conditions \reff{eq:2*k}
   with respect to $Y^{k+1}$ can be written as follows:
 \begin{eqnarray*} 
 	Y^{k+1}_{j}(x_j,t)= \left[RZ^{k+1}\right]_j(t)-\left[RZ^{k}\right]_j(t)
 	+\chi_j^k(x_j,t), \;\;\;  j\le n,
 \end{eqnarray*}  
 or, in the above notation, as
\beq\label{bcY}
Y^{k+1}_{j}(x_j,t)= \left[R\left(\bar Y^{k+1}\right)\right]_j(t)+\zeta^k_j(t), \;\;\; j\le n.
\ee 
Therefore,  the function $w^{k+1}$ is the classical $BC^2$-solution to the  system
 \begin{eqnarray}
 & & \d_ tw^{k+1} + A(x,t,V^k) \partial_ x w^{k+1}+ B(x,t,V^k) w^{k+1} = f^k(x,t)
 \label{th-5}
 \end{eqnarray}
 with the boundary conditions \reff{zam2}, \reff{bcY}, where 
 \begin{eqnarray*}
  f^k(x,t)   &=& - \left(B^k -   B^{k-1}\right)V^k- \left(A^k - 
   A^{k-1}\right)\partial_ xV^k\\
 &=& - \int_0^1\partial_3 B\left(x,t,\sigma V^k(x,t)+(1-\sigma)V^{k-1}(x,t)\right)d\sigma
w^k(x,t)V^k(x,t)\\ &&-\int_0^1\partial_3 A\left(x,t,\sigma V^k(x,t)+(1-\sigma)V^{k-1}(x,t)\right)d\sigma
w^k(x,t)\partial_ xV^k(x,t).
\end{eqnarray*}

Now we show that the sequence $w^{k+1}$ converges to zero in $BC^1_t(\Pi;\R^n).$
By Claim \ref{ots5}, the functions  $V^k$ and  $V^{k-1}$ satisfy  the same estimate \reff{apr_seq}. 
 On the account of \reff{phiest}--\reff{minus}, there exists a 
 constant $N_1$  not depending on $V^k$, $V^{k-1}$, and   $w^k$ such that
  \beq\label{app1}
  \begin{array}{rcl}
  \|f^{k}\|_{BC^1_t} &\le& N_1 \left(\|V^k\|_{BC^1} + \|\partial_{x} V^k\|_{BC^1_t}\right)\|w^k\|_{BC^1_t}\\
&  \le&  N_1 K_2 \left(\| f\|_{BC_t^2} + \|h\|_{BC^2}\right) \|w^k\|_{BC^1_t}.
  \end{array}
\ee
Similarly we obtain  the bound
  \begin{eqnarray}\label{app2}
   \|\zeta^k\|_{BC^1} \le  N_1 K_2 \left(\| f\|_{BC_t^2} + \|h\|_{BC^2}\right) \|w^k\|_{BC^1_t},
  \end{eqnarray} 
  where the constant $N_1$  does not depend on $V^k$, $V^{k-1}$, and   $w^k$
  and is chosen to satisfy both the inequalities  \reff{app1}  and  \reff{app2}. 
   By Part 1 of Theorem~\ref{prop3}, the solution   $w^{k+1}$
  to the problem \reff{th-5}, \reff{zam2}, \reff{bcY} satisfies the estimate
  \reff{ots55} with $v$, $f$, and $h$ replaced by $w^{k+1}$, $f^k$, and
  $\zeta^k$, respectively. Combining this estimate with
  \reff{app1}--\reff{app2}, we derive  the
inequality
\beq\label{fin1}
 \begin{array}{rcl}
\|w^{k+1}\|_{BC^1_t}& \le& K_1\left(\|f^k\|_{BC^1_t }+\|\zeta^k\|_{BC^1}\right)\\
&\le&
\displaystyle K_1K_2N_1 \left(\|f\|_{BC^2_t}+\|h\|_{BC^2}\right) \|w^k\|_{BC^1_t}.
 \end{array}
\ee
Set
\beq\label{eps1}
\eps_2=\min\left\{\delta_1/ K_2,(K_1K_2N_1)^{-1}\right\}.
\ee
If
 $ \|f\|_{BC^2_t}+\|h\|_{BC^2} < \eps_2,$
then,  due to \reff{fin1}, the sequence $w^k$ is strictly contracting  in  $BC^1_t(\Pi;\R^n)$ and, hence  tends to zero in  $BC^1_t(\Pi;\R^n).$

By the inequality \reff{apr_seq} and the assumptions of Theorem~\ref{prop3},
the inverse $(A^k)^{-1}$ exists for every $k$ and, moreover, is bounded in
 $BC(\Pi;\R^n)$ uniformly in $k$.
Now, the equation (\ref{th-5}) yields
\beq\label{parx}
\begin{array}{rcl}
 && \|\partial_x w^{k+1}\|_{BC}\le \|(A^k)^{-1}\|_{BC}\left(
 \|f^k\|_{BC} +  \|\partial_t w^{k+1}\|_{BC} +
  \|B^{k}\|_{BC}\|w^{k+1}\|_{BC}\right) \nonumber\\ [2mm]
 &&  \le K_2 N_1\|(A^k)^{-1}\|_{BC}\left(1 + K_1 + K_2 \|B^{k}\|_{BC}\right) 
 \left(\|f\|_{BC^2_t} + \|h\|_{BC^2}\right) \|w^k\|_{BC^1_t},
\end{array}
\ee
which together with \reff{fin1} gives the convergence $\|w^{k+1}\|_{BC^1} \to 0$ as $k \to \infty$.

Finally, because of \reff{zam2*},
the sequence $V^k$ converges to some function $V^*$ in $BC^1(\Pi;\R^n)$.
It is a simple matter to show that  $V^*$ is a classical solution to the problem 
(\ref{eq:1})--\reff{z*}.
The proof of the claim is complete.
\end{subproof}

\begin{cclaim} \label{convC2}
There exist positive constants $\varepsilon\le\eps_2$ and $\delta \le\de_1$ such that, if $\|f\|_{BC^2_t}+\|h\|_{BC^2} \le \varepsilon,$
then
 the classical solution $V^*$ belongs to $BC^2(\Pi;\R^n)$ and satisfies the estimate
 \begin{eqnarray}
 \| V^*\|_{BC^2_t}+ \|\partial_{x} V^*\|_{BC^1_t}\le \delta.
\label{*1**}
\end{eqnarray}
\end{cclaim}
\begin{subproof}
we start with proving that the sequence $V^k$ converges in $BC^2(\Pi;\R^n)$.
First show that the sequence 
$W^{k+1} = (Q^k)^{-1} \d_tV^{k+1}=\d_tU^{k+1}+
(Q^k)^{-1}\left(\d_tQ^k + \d_3Q^kQ^{k-1}W^k\right)U^{k+1}$
converges in $BC^1_t(\Pi;\R^n)$. To this end, we
differentiate the problem \reff{th-2}--\reff{zamk}
with respect to $t$ and, similarly to \reff{eq:1120} and \reff{eq:1155}, write down the resulting problem
in the diagonal form with respect to $W^{k+1}$, as follows:
\beq\label{eq:1120*}
\begin{array}{rr}
 \partial_t W^{k+1} + \hat a(x,t,V^k)\partial_x W^{k+1} 
+ b^{1}(x,t,V^k)W^{k+1}\\ [2mm]=g^{1k}(x,t,W^k)W^{k+1}+g^{2k}(x,t,W^k),
 \end{array}
\ee
\begin{equation} \label{eq:1155*}
 W^{k+1}_{j}(x_j,t)
  = (\widetilde Ry^{k+1})_j(t) + h^{k}_j(t,W^k), \  j\le n,
\end{equation}
 where
 \begin{eqnarray*} 
  b^{1}(x,t,V^k) &=&  b^k - (Q^k)^{-1} \d_t A^k (A^k)^{-1}Q^k, 
 \\ [1mm]
  g^{1k}(x,t,W^k) &=& (Q^k)^{-1}\d_3 A^k Q^{k-1} W^k (A^k)^{-1}Q^k - (Q^k)^{-1} \d_3 Q^k Q^{k-1} W^k  \\
  &&- (Q^k)^{-1}A^k \d_3 Q^k (A^{k-1})^{-1}\left( f - Q^{k-1}W^k - B^{k-1}V^k\right),
 \\ [1mm]
 g^{2k}(x,t,W^k) &=& (Q^k)^{-1}\Bigl(-\d_tB^kV^{k+1}+\d_tA^k\,(A^k)^{-1} B^kV^{k+1}
 \\ &&+ \d_t f - \d_tA^k \,(A^k)^{-1}f  - \d_3 B^k Q^{k-1} W^kV^{k+1} \\ 
 &&+ \d_3 A^k\,Q^{k-1} W^k (A^k)^{-1}(B^k V^{k+1} -f) \Bigl), \\[1mm]
h^{k}_j(t,W^k) & = & \left(R^\prime Z^{k+1}\right)_j(t) - (\widetilde R\rho^{k})_j(t)+ h^\prime_j (t)
+ \varrho_j^{k}(x_j,t),   \   j\le n,  \\
   y^{k+1}(t)&=&\left(W^{k+1}_1(1,t),\dots,W^{k+1}_{m}(1,t),W^{k+1}_{m+1}(0,t),\dots,W^{k+1}_{n}(0,t)\right),\\
   \rho^{k}(t)&=&\Bigl(\varrho_1^{k}(1,t),\dots,\varrho_m^{k}(1,t), 
            \varrho_{m+1}^{k}(0,t),\dots, \varrho_{ n}^{k}(0,t)\Biggr), \\
             \varrho_j^{k}(x,t) & = & \left[(Q^k)^{-1}(\d_t Q^k+ \d_3Q^kQ^{k-1}W^k)U^{k+1} \right]_j(x,t).
 \end{eqnarray*} 
  It is evident that the sequence 
$W^{k+1}$ of solutions to the problem \reff{eq:1120*}--\reff{eq:1155*}
converges in $BC^1_t(\Pi;\R^n)$ if and only if the sequence $W_t^{k+1}=\d_tW^{k+1}$
converges in $BC(\Pi;\R^n)$.  To prove  the last statement, 
we differentiate the system \reff{eq:1120*} in the distributional sense and the boundary 
 conditions \reff{eq:1155*} pointwise in $t$. We, therefore,  obtain the following problem with respect to
 $W_t^{k+1}$:
\begin{eqnarray}
  \partial_t W_t^{k+1} + a^k\partial_x W_t^{k+1} +b^{2k}W_t^{k+1}
={\cal G}_{1}(k)W_t^{k+1}+{\cal G}_{2}(k)W_t^{k}+g^{3k}, \label{eq:1120**} \\[1mm]
 (W_t)^{k+1}_{j}(x_j,t)
  = (\widehat R(y^{k+1})^\prime)_j(t)+ \left[{\cal H}(k)W_t^{k}\right]_j(t)+\tilde h^{k}_j(t), \ j\le n,  \label{eq:1155**}
\end{eqnarray}
 where
  \begin{eqnarray*} 
 	b^{2k}&=&b^{1k}-\d_ta^k\,(a^k)^{-1},\\ 
 	g^{3k}&=&\Bigl( \d_tg^{1k} - \d_3 b^{1k}Q^{k-1}W^k  \\ &&+ 
 	\left(\d_3 a^{k}Q^{k-1}W^k (a^k)^{-1} + \d_t a^{k} (a^k)^{-1}\right) (b^{1k} -  g^{1k}) - \d_t b^{1k} \Bigl) W^{k+1}  \\
 	&&	
 	+ \d_tg^{2k} -  \d_3 a^{k}Q^{k-1}W^k (a^k)^{-1} g^{2k} - \d_t a^{k} (a^k)^{-1} g^{2k} , \\ [2mm]
 	\tilde h^{k}_j(t)&=&(\widetilde R^\prime y^{k+1})_j(t) + 
 	\d_th^{k}_j(t,W^k).
 \end{eqnarray*}  
  Moreover, $b^{1k}$ is used to denote the function $b^{1}(x,t,V^k)$, while
    the operators
 ${\cal G}_{1}(k),{\cal G}_{2}(k)\in\LL(BC(\Pi;\R^n))$ and
  ${\cal H}_j(k)\in\LL(BC(\Pi;\R^n);BC(\R;\R^n))$ are defined by
 \begin{eqnarray*} 
  \left[{\cal G}_{1}(k)W_t^{k+1}\right](x,t) &=&\left(g^{1k} + 
 \partial_3 a^k Q^{k-1}W^k(a^k)^{-1}\right) W_t^{k+1},\\ 
  \left[{\cal G}_{2}(k)W_t^{k}\right](x,t)&=&\d_3g^{1k}W_t^{k}W^{k+1}+\d_3g^{2k}W_t^{k},\\ 
 \left[{\cal H}(k)W_t^{k}\right]_j(t)&=&\d_2h^{k}_j(t,W^k)W_t^k(x_j,t),\ \ j\le n.
 \end{eqnarray*} 

 Similarly to Claim \ref{cl5} in the proof of Theorem \ref{smooth_sol},
 the function $W_t^{k+1}$ satisfies \reff{eq:1120**} in the distributional sense
 and \reff{eq:1155**} pointwise if and only if it satisfies the following
 operator equation:
  \beq \label{oper3*}
  \begin{array}{rcl}
  W_t^{k+1}& =& C(k)W_t^{k+1} + D(k) W_t^{k+1}  \\
  &&	+ 
  F(k)\left({\cal G}_{1}(k)W_t^{k+1}+{\cal G}_{2}(k)W_t^{k}+g^{3k},{\cal H}(k)W_t^{k}+\tilde h^{k}\right),
   \end{array}
 \ee
where the operators
 $C(k)$, $D(k)$, and $F(k)$ are defined by the right hand sides of
 the corresponding formulas in
\reff{CDF}  with $a$, $b$, 
 and $R$ replaced, respectively,  by
$a^k$, $b^{2k}$, and $\widehat R$.
Moreover, the functions  $\om_j$, $c_j$, and $d_j$ are replaced appropriately by 
$\om_j^k$, $c_j^k$, and $d_j^k$. Note that computing $C(k) W^{k+1}_t$,
we put $z=(y^{k+1})^\prime$ in the right hand side of the first formula 
in~\reff{CDF}.

 Iterating (\ref{oper3*}), we get
 \beq \label{oper4*}
 \begin{array}{rcl} 
 && W^{k+1}_t   =  C(k)W^{k+1}_t + \left(D(k) C(k) + D^2(k)\right)W^{k+1}_t \nonumber\\ 
   &&\ \ +(I+D(k)) F(k)\left({\cal G}_{1}(k)W_t^{k+1}+{\cal G}_{2}(k)W_t^{k}+g^{3k},{\cal H}(k)W_t^{k}+\tilde h^{k}\right).
 \end{array} 
 \ee
 Now we intend to show that there exists $\de_2\le\de_1$ such that,
 given a nonnegative integer $k$ and
 $V^k$ satisfying the estimate \reff{apr_seq} with $\de_2$
 in place of $\de_1$, the  formula \reff{oper4*} is equivalent
 to the following one:
 \beq \label{oper5*}
 W^{k+1}_t  = {\cal A}(k)W^{k}_t+X^k,
 \ee
 where ${\cal A}(k)\in\LL(BC(\Pi;\R^n))$ and  $X^k\in BC(\Pi;\R^n)$
 are given by
 \beq \label{oper6*}
 \begin{array}{rcl}
  {\cal A}(k)W  &=& 
  \left[I-C(k)-(I+D(k)) F(k)\left({\cal G}_{1}(k),0\right)\right]^{-1}\\ [2mm]
  &&\times(I+D(k)) F(k)
   \left({\cal G}_{2}(k)W,{\cal H}(k)W\right),\\ [2mm]
  X^k&=&\left[I-C(k)-(I+D(k)) F(k)\left({\cal G}_{1}(k),0\right)\right]^{-1}
  \left(D(k) C(k) + D^2(k)\right)W^{k+1}_t\\ [1mm]
   &+&\left[I-C(k)-(I+D(k)) F(k)\left({\cal G}_{1}(k),0\right)\right]^{-1}(I+D(k)) F(k)
   \left(g^{3k},\tilde h^{k}\right).
 \end{array} 
 \ee
 It suffices to show that, for every
  nonnegative integer $k$,  the operator $I-C(k)-(I+D(k)) F(k)\left({\cal G}_{1}(k),0\right)$
 is invertible and has a bounded inverse. 
  Even more, we will show
 that the inverse is bounded uniformly in $k$. With this aim, denote by
$G_0(k)$   operator defined by the right hand sides of  \reff{Ctilde}, where $a_j$, $b_{jj}$, and $\om_j$ are replaced, respectively,
 by $a_j^k$, $b_{jj}^{2k}$, and $\om_j^k$. 
Moreover, denote by $G_2(k)$  operator defined by the right hand sides of  the second formula in \reff{Ctilde1}, where $a_j$, $b_{jj}$, and $\om_j$ are replaced, respectively,
 by $a_j^k$, $b_{jj}^{k}$, and $\om_j^k$. 
 
 Note that, similarly to \reff{bjj}, we have $
 b_{jj}^{1k}=b_{jj}^{k}-(a_{j}^{k})^{-1}\d_ta_{j}^{k}.
 $
Then, accordingly to the notation introduced above, the function $b_{jj}^{2k}$
 is given by the formula
 $
 b_{jj}^{2k}=b_{jj}^{k}-2(a_{j}^{k})^{-1}\d_ta_{j}^{k}.
 $
This means that the operators $G_0(k)$ and $G_2(k)$ coincide.

Therefore, on the account of the  estimates \reff{apr_seq} and  \reff{G2},  the operators
$G_2(k)$ fulfill the inequality  $\|G_2(k)\|_{\LL(BC(\R,\R^n))}<1-\nu_3$
for all $k\in\N$ and, hence,  the inequality  $\|G_0(k)\|_{\LL(BC(\R,\R^n))}<1-\nu_3$
for all $k\in\N$.

Finally, similarly to the proof of
  the invertibility  of   $I-C$ in Subsection~\ref{bounded_sol2}, 
   the invertibility   of  $I-C(k)$  follows from the invertibility of 
    $I-G_0(k)$  (see the inequality \reff{uj}).
   Furthermore, the following estimate is true for all $k\in\N$:
 $$
\|(I-C(k))^{-1}\|_{\LL(BC(\Pi;\R^n))}\le 1+\nu_3^{-1}\|C(k)\|_{\LL(BC(\Pi;\R^n))}.
$$
As the operators $C(k)$ are  bounded uniformly in $k$, the inverse operators
$(I-C(k))^{-1}$ are  bounded uniformly in $k$ also.
Taking into account that the set of all invertible 
operators whose inverses are bounded  is open, our task is, therefore, reduced to show that
the operator $(I+D(k)) F(k)\left({\cal G}_{1}(k),0\right)$ is  sufficiently small 
 whenever $\de_1$ is sufficiently small. Note that  Claim \ref{ots5}
 is true with $\de_2$
 in place of $\de_1$
  for any $\de_2\le\de_1$.  This implies that  for  any $\sigma>0$ there is 
$\de_2$ such that for all $V^k$ fulfilling
\reff{apr_seq}  with $\de_2$
in place of $\de_1$, we have
$\|{\cal G}_1(k)\|_{\LL(BC(\Pi;\R^n))}=\|g^1(x,t, W^k(x,t))\|_{BC}\le\sigma$ for all $k$. 
Moreover,  the operators $D(k)$ and  
$F(k)$ are bounded, uniformly in $k$.
Consequently, if $\de_2$ is sufficiently small, then for all $k\in\N$ and all $f$ and $h$ satisfying \reff{sm7aa} with $\de_2$ in place of $\de_1$, the operator 
$I-C(k)-(I+D(k)) F(k)\left({\cal G}_{1}(k),0\right)$
is invertible and the inverse is  bounded by a constant not depending on $k$.
Fix  $\de_2$ satisfying the last property.
The equivalence of \reff{oper4*} and \reff{oper5*} is, therefore, proved.
 
Now, to prove that the sequence $W_t^{k+1}$
converges in $BC(\Pi;\R^n)$ as $k\to\infty$,  we apply  to  the equation \reff{oper5*} a linear version of the 
fiber contraction principle, see  \cite[Lemma A.1]{Hopf}.
Accordingly to \cite[Lemma A.1]{Hopf}, we have to show that, first,
\beq\label{oper7*}
X^k \mbox{ converges in } BC(\Pi;\R^n) \mbox{ as } k\to\infty,
\ee
second, that there exists $c<1$  such that for all $W\in BC(\Pi;\R^n)$ it holds 
\beq\label{oper8*}
\|{\cal A}(k)W\|_{BC}\le c\|W\|_{BC},
\ee
and, third, that
\beq\label{oper9*}
{\cal A}(k)W \mbox{ converges in } BC(\Pi;\R^n) \mbox{ as } k\to\infty.
\ee

To show \reff{oper7*}, note that the operators 
$D(k)$ and $C(k)$  depend  neither on $W^k_t$ nor on 
$W^k$. Similarly to the proof of Claim \ref{2} in Section \ref{sec:regul}, one can show that the operators $D(k) C(k)$ and $D(k)^2$ are smoothing and map $BC(\Pi;\R^n)$ into $BC^1_t(\Pi;\R^n)$.
This implies that $D(k) C(k)W_t^{k+1}$ and $D(k)^2W_t^{k+1}$, 
actually, do not depend on $W_t^{k+1}$,
but on $W^{k+1}$. Moreover, using  \reff{apr_seq},
we get the following  estimate:
\begin{eqnarray*} 
\left\|\left(D(k)C(k) + D(k)^2\right)W_t^{k+1}\right\|_{BC} \le \widehat K_{11} \|W^{k+1}\|_{BC}
\end{eqnarray*}
for some $\widehat K_{11}$ not depending on $k$.
It follows  that the right hand side of the second formula in \reff{oper6*} does not depend on $W_t^{k}$ for all $k$ and, therefore, the convergence \reff{oper7*} immediately follows
from Claim \ref{convV}. 

Since all the operators in the right hand side of the first formula in \reff{oper6*} do not depend on $W_t^{k}$ for all $k$,
the convergence \reff{oper9*} follows again from Claim \ref{convV}.

It remains to prove \reff{oper8*}.  Due to Claim \ref{ots5}, for any $\sigma>0$
there exists $\delta\le\de_2$ such that for all $f$ and $h$ satisfying the bound \reff{sm7aa} 
 with $\de$ in place of $\de_1$ (and, hence for
$V^k$ satisfying the bound \reff{apr_seq}  with $\de$
in place of $\de_1$) it holds
$\|{\cal G}_2(k)\|_{\LL(BC(\Pi;\R^n))}+\|{\cal H}(k)\|_{\LL(BC(\Pi;\R^n))}\le\sigma$.
Moreover,  if $\de\le\de_2$ is sufficiently small, then
all other operators in the right hand side of 
the first equality in \reff{oper6*} are bounded uniformly in $k$. 
We, therefore,  conclude  that there exists $\de\le\de_2$ such
that \reff{oper8*} is fulfilled. 

Set 
$
\eps=\min\left\{\delta/ K_2,(K_1K_2N_1)^{-1}\right\}
$
(see also \reff{eps1}).
Therefore, by Lemma \cite[Lemma A.1]{Hopf}, if $ \|f\|_{BC^2_t}+\|h\|_{BC^2} < \eps,$ then
the  sequence $W_t^{k+1}$
converges in $BC(\Pi;\R^n)$ as $k\to\infty$.

Finally, by  Claim \ref{convV} and the equality 
$W^{k+1} = (Q^k)^{-1} \d_tV^{k+1}$, we conclude that the second derivative 
of $V^*$ in $t$ exists and that the sequence $\d^2_tV^k$ converges to $\d^2_tV^*$
in $BC(\Pi;\R^n)$
as $k\to\infty$.
Differentiating \reff{th-2} first in $t$ and then in $x$
and using the fact that $A(x,t,V^k)^{-1}$  is bounded uniformly in $k$, we conclude that $V^k$ converges to $V^*$  in $BC^2(\Pi;\R^n)$ as $k\to\infty$. 

The desired 
estimate \reff{*1**} now easily follows from the bound \reff{apr_seq} and the system \reff{th-2} differentiated in $t$.
 \end{subproof}
 
 \begin{cclaim} 
	Let $\eps$ and $\de$ be as in Claim \ref{convC2}.
	Then there exists  $\delta^\prime=\delta^\prime(\eps,\delta)$ such that,  if  $\|f\|_{BC_t^2} +\|\d_x f\|_{BC} + \|h\|_{BC^2}\le\eps$, then  
	$\|V^*\|_{BC^2}\le~\delta^\prime$.
\end{cclaim}

\begin{subproof}
	Consider  the system \reff{eq:1} with $V$ replaced by $V^*$,
differentiated  in~$x$. Using the fact
that $A(x,t,V^*)$  has a bounded inverse and taking into account
 the estimate \reff{*1**},
 we derive an
upper bound for
$\d_x^2 V^*$. One can easily see that  this bound depends  on $\eps$ and $\delta$
and that  there 
exists   $\delta^\prime=\delta^\prime(\eps,\delta)$ such that
$\|V^*\|_{BC^2}\le \delta^\prime$, as desired.	
 \end{subproof}

\begin{cclaim} Let $\eps$ and $\delta$ be as in Claim \ref{convC2}. Then for any $f$ and $h$ such that $\|f\|_{BC^2_t}+\|h\|_{BC^2} \le\eps$,
 the classical solution to the problem  (\ref{eq:1})--\reff{z*} 
 fulfilling the estimate (\ref{*1**}) is unique.
\end{cclaim}

\begin{subproof}
On the contrary, suppose  that $\widetilde V$ is a classical  solution to the problem 
(\ref{eq:1})--\reff{z*} different from $V^*$,
such that $ \|\widetilde V\|_{BC^2_t} + \|\d_x\widetilde V\|_{BC^1_t} \le \delta.$ Then, due to \reff{phiest}--\reff{minus}, the functions $\widetilde A(x,t) = A(x,t,\tilde V(x,t))$, $\widetilde Q(x,t) = Q(x,t,\tilde V(x,t))$
and $\widetilde B(x,t) = B(x,t,\tilde V(x,t))$
fulfill the inequalities
$$
\begin{array}{cc}
\|\widetilde A - A^0\|_{BC^2_t} + \|\d_x \widetilde A- \d_x A^0\|_{BC^1_t}\le \varepsilon_1, \quad 
\|\widetilde B - B^0\|_{BC^1_t} \le \varepsilon_1\\ [2mm]
\|\widetilde Q - Q^0\|_{BC^2_t} + \|\d_x \widetilde Q- \d_x Q^0\|_{BC^1_t}\le \varepsilon_1.\end{array}$$
The difference $\tilde w^{k+1} =V^{k+1}- \widetilde V$ satisfies the  system
$$
\partial_t \tilde w^{k+1}  + \widetilde A(x,t) \partial_x \tilde w^{k+1}  + \widetilde B(x,t) \tilde w^{k+1}  = \tilde f^{k}(x,t)
$$
 and the boundary conditions \reff{zam2}, (\ref{bcY})
 with $w^{k+1}$ replaced  by  $\tilde w^{k+1}$ and with
\begin{eqnarray*}
Y^{k+1}= \widetilde Q^{-1}\tilde w^{k+1},\quad 
\chi^k=\widetilde Q^{-1}(Q^k - \widetilde Q)  U^{k+1}.
\end{eqnarray*}
Here  $\widetilde Q(x,t) =Q\left(x,t,\widetilde V(x,t)\right)$
and
  \begin{eqnarray*}
  \tilde f^{k}(x,t)
  = \left(\widetilde B(x,t)-B^{k}(x,t) \right)V^{k+1}(x,t) +
 \left(\widetilde A(x,t)-A^{k}(x,t)\right)\partial_x V^{k+1}(x,t).
 \end{eqnarray*}
By the argument as in the proof of Claim \ref{convV},  the functions $\tilde f^{k}(x,t)$ and $\zeta^{k}(x,t)$
are $C^1$-smooth in $t$ and
satisfy the  upper bounds \reff{app1} and \reff{app2}
with $f^k$ and $w^k$ replaced by $\tilde f^k$ and $\tilde w^k$,
respectively.

Similarly to \reff{fin1} and \reff{parx}, we derive the bounds
\begin{eqnarray*}
& & \|\tilde w^{k+1}\|_{BC^1_t}   \le
 K_1 K_2 N_1 \left(\|f\|_{BC^2_t}+\|h\|_{BC^2}\right) \|\tilde w^k\|_{BC^1_t},\\
&& \|\partial_x \tilde w^{k+1}\|_{BC} \displaystyle
 \le \frac{1}{\Lambda_0}K_2 N_1\left(1 + K_1 + K_2 \|B^{k}\|_{BC}\right) 
 \left(\|f\|_{BC^2_t} +\|h\|_{BC^2}\right) \|\tilde w^k\|_{BC^1_t}.
\end{eqnarray*} 
The desired  convergence $\|\tilde w^{k}(t)\|_{BC^1} \to 0$
as $k \to \infty$ follows. This means that $\widetilde V(x,t) = V^*(x,t)$, contradicting to our assumption.
\end{subproof}

\subsection{Proof of Part 2 of Theorem \ref{main}: Almost periodic solutions}\label{ap}
We have to prove that the constructed solution $V^*(x,t)$
 is  Bohr almost periodic in $t$.
The proof uses the fact that the limit of a uniformly convergent sequence of Bohr almost periodic functions 
 depending uniformly on parameters
 is almost periodic uniformly in parameters \cite[p. 57]{Cord}.
Moreover,  we will use the fact that,  if a function $w(x,t)$ has bounded and continuous partial derivatives up to the second order in both
 $x \in [0,1]$ and in $t \in \mathbb R$ and  is Bohr almost periodic in $t$ uniformly in $x$ (or, simply,  almost periodic), the last property is true for
 $\partial_x w(x,t)$
 and $\partial_t w(x,t)$ also. Specifically, the almost periodicity of
 $\partial_t w(x,t)$ follows from \cite[Theorem 2.5]{Cord}, while the almost periodicity of
 $\partial_x w(x,t)$ is shown in \cite[Section 5.2]{KRT2}.
We are, therefore,  reduced  to showing that the approximating sequence $V^k$,  constructed in Section~\ref{sec:bounded},
is a sequence of    almost periodic functions.

We  use the induction on $k$.
Recall that $V^0\equiv 0$.
Assuming that the iteration $V^k(x,t)$
is Bohr almost periodic for an arbitrary fixed $k\in\N$, 
let us prove that $V^{k+1}(x,t)$ is almost periodic also. By the assumptions of  the theorem, 
the matrices $ A(x,t,V^k(x,t)),$ 
$ B(x,t,V^k(x,t)),$  $Q(x,t,V^k(x,t))$,
$\d_xQ(x,t,V^k(x,t))$, and $\d_tQ(x,t,V^k(x,t))$ are almost periodic  as compositions of almost periodic functions. 
Below we will use a slightly modified notation for $\hat a(x,t,V^k)$  and $\hat b(x,t,V^k)$ (see \reff{hat}), namely
$a^k(x,t) = \hat a(x,t,V^k(x,t))$ and 
$b^k(x,t) = \hat b(x,t,V^k(x,t))$. Set 
$q^k(x,t) = Q(x,t,V^k(x,t))$.
It follows that $a^k$ and $b^k$ are almost periodic.
Fix  $\mu > 0$ and let $\tau$ be a $\mu$-almost period of the  matrices 
$a^k, q^k$ and $b^k$.
Then the differences $\tilde a^k(x,t) = a^k(x,t + \tau) - a^k(x,t),$ 
$\tilde b^k(x,t) = b^k(x,t + \tau) - b^k(x,t)$, and $\tilde q^k(x,t) = q^k(x,t + \tau) - q^k(x,t)$
satisfy the inequalities
\begin{eqnarray} \label{tilde}
 \|\tilde a^k\|_{BC} \le \mu,\quad \|\tilde b^k\|_{BC} \le \mu,\quad \|\tilde q^k\|_{BC} \le \mu
\end{eqnarray}
uniformly in $x$ and $t$. 

First derive a few   simple estimates. Let $\om_j^k(\xi,x,t)$ be the solution to the equation
 \reff{char} where $a_j$ is replaced by $a_j^k$.
Then the following  identity is true:
 \begin{eqnarray*}
  \frac{d}{d\eta}\left( \om_j^k(\eta,x,t) + \tau -  \om_j^k(\eta,x,t+\tau)\right) = \frac{1}{a_j^k(\eta,\om_j^k(\eta,x,t))} -
\frac{1}{a_j^k(\eta,\om_j^k(\eta,x,t+\tau))}.
 \end{eqnarray*}
Since $\om_j^k(x,x,t) = t$ and $\om_j^k(x,x,t+\tau) = t+\tau,$ it holds
 \beq\label{id}
  \begin{array}{ll}
 	  \displaystyle \om_j^k(\eta,x,t) + \tau -  \om_j^k(\eta,x,t+\tau)\nonumber \\ \ \  =
 \displaystyle \  \int_x^\eta \Biggl(\frac{1}{a_j^k(\xi,\om_j^k(\xi,x,t))} - \frac{1}{a_j^k(\xi,\om_j^k(\xi,x,t+\tau))}\Biggr)d\xi \nonumber \\ [5mm]
    \displaystyle \ \ = \int_x^\eta \frac{a_j^k(\xi,\om_j^k(\xi,x,t+\tau)) - a_j^k(\xi,\om_j^k(\xi,x,t)+\tau)}{a_j^k(\xi,\om_j^k(\xi,x,t))
a_j^k(\xi,\om_j^k(\xi,x,t+\tau))}d\xi  \nonumber\\ [5mm]
   \displaystyle \ \  \ \ \ \ + \int_x^\eta \frac{a_j^k(\xi,\om_j^k(\xi,x,t)+\tau) - a_j^k(\xi,\om_j^k(\xi,x,t))}{a_j^k(\xi,\om_j^k(\xi,x,t))
a_j^k(\xi,\om_j^k(\xi,x,t+\tau))}d\xi. 
 \end{array}\ee
 By \reff{tilde},
 $$
 |a_j^k(\xi,\om_j^k(\xi,x,t)+\tau) - a_j^k(\xi,\om_j^k(\xi,x,t))|\le\mu,
 $$
 the estimate being uniform in $\xi,x,t$, and $j$. 
 Due to the mean value theorem,
 $$ \begin{array}{ll }
 a_j^k(\xi,\om_j^k(\xi,x,t+\tau)) - a_j^k(\xi,\om_j^k(\xi,x,t)+\tau)
 =\displaystyle
 ( \om_j^k(\xi,x,t+\tau)-\om_j^k(\xi,x,t) - \tau)\\ \qquad\qquad\times \displaystyle \int_0^1\d_2a_j^k\left(\xi,\al\om_j^k(\xi,x,t+\tau)+(1-\al)(\om_j^k(\xi,x,t)+\tau
  \right)d\al.
\end{array} 
 $$
 Applying
  the Gronwall's inequality to the  identity \reff{id}, we derive the estimate
 \begin{eqnarray} \label{ap1}
  \left|\om_j^k(\eta,x,t) + \tau - \om_j^k(\eta,x,t+\tau)\right| \le \frac{\mu}{\Lambda_0^2}
\exp\Biggl\{\frac{\|a_j^k\|_{BC^1_t}}{\Lambda_0^2}\Biggr\}= L_1\mu,
 \end{eqnarray}
the constant $L_1$ being independent of $\mu$, $\eta,x,t$, and $j$.

 Next we  show that $a_j^k(\eta, \om_j^k(\eta,x,t))$ and $b_{ji}^k(\eta, \om_j^k(\eta,x,t))$
are almost periodic. For that we use
(\ref{ap1}) and the fact that  $\tau$ is a $\mu$-almost period  of $a_j^k$ and  $b_{ji}^k$. We get
\beq\label{ap3}
 \begin{array}{ll}
 \left|a_j^k(\eta, \om_j^k(\eta,x,t)) - a_j^k(\eta, \om_j^k(\eta,x,t+\tau))\right| \\
 [2mm] \qquad\le
  \left|a_j^k(\eta, \om_j^k(\eta,x,t)) - a_j^k(\eta, \om_j^k(\eta,x,t) + \tau)\right| \\ [2mm]\qquad\quad
 + \left|a_j^k(\eta, \om_j^k(\eta,x,t)+\tau) - a_j^k(\eta, \om_j^k(\eta,x,t+\tau))\right| \\ [2mm] \qquad\le
 \left(1 + L_1\|\partial_t a_j^k\|_{BC}\right)\mu \le L_2 \mu,
 \end{array}
\ee
where $L_2$ does not depend on $\mu$,  $\eta,x,t$, and $j$.
 Similar estimates are true for $b_{ji}^k$ and $q_{ji}^k$, namely
 \beq\label{ap4}
 \begin{array}{ll}
   |b_{ji}^k(\eta, \om_j^k(\eta,x,t)) - b_{ji}^k(\eta, \om_j^k(\eta,x,t+\tau))| \le L_2 \mu,\\ [2mm]
    |q_{ji}^k(\eta, \om_j^k(\eta,x,t)) - q_{ji}^k(\eta, \om_j^k(\eta,x,t+\tau))| \le L_2 \mu,
    \end{array}
\ee
where  $L_2$  is chosen to be a  common constant satisfying both 
\reff{ap3} and \reff{ap4}.

 Now we prove that
 \beq\label{ap2}
(Rv)_j(\om_j^k(x_j,x,t)) \in AP(\Pi)\quad\mbox{for all}\quad
v\in AP(\R, \mathbb R^n)\cap BC^1(\R,\R^n).
\ee
  Fix an arbitrary 
$v\in AP(\R, \mathbb R^n)\cap BC^1(\R,\R^n)$.
By  the assumption, $(Rv)(t)\in AP(\R, \mathbb R^n)$.
Let $\tau$ be a common   $\mu$-almost period  of the  functions 
$(Rv)(t)$ and $a^k(x,t)$. 
Applying the mean value theorem and using the assumption $\bf{(A3)}$ and the estimate \reff{ap1},
we get
$$
\begin{array}{ll}
\left|(Rv)_j(\om_j^k(x_j,x,t)) - (Rv)_j(\om_j^k(x_j,x,t+\tau))  \right|
\\ [2mm]
\quad\quad\le
\left|(Rv)_j(\om_j^k(x_j,x,t)) - (Rv)_j(\om_j^k(x_j,x,t)+\tau)  \right|\\ [2mm]
\quad\quad\quad+ \left|(Rv)_j(\om_j^k(x_j,x,t)+\tau)-(Rv)_j(\om_j^k(x_j,x,t+\tau))  \right|
\\ [2mm]
\quad\quad\le \displaystyle\mu\left(1 + L_1
\sup_{t\in\R}\left|\frac{d}{dt}(Rv)_j(t)\right|
\right),
\end{array}
 $$
 which proves \reff{ap2}.

 The estimates \reff{ap3} and \reff{ap4} imply that the
 functions in the right hand sides of the equalities in
 \reff{c0} 
  with  $a_j^k$, $b_{jj}^k$, and $\om_j^k$ in place of $a_j$,
$b_{jj}$, and $\om_j$, respectively,
are almost periodic for all $j \le n$, uniformly in $\xi,x \in [0,1]$.
Let $\widehat C(k)$, $\widehat D(k)$, and $\widehat F(k)$ be defined by the right hand side of
\reff{CDF}  with  $a_j$, $b_{jj}$, and $\om_j$ replaced by  $a_j^k$,
$b_{jj}^k$, and $\om_j^k$, respectively.
Taking into account \reff{ap2}, we conclude that the operators $\widehat C(k)$, $\widehat D(k)$, and $\widehat F(k)$ map the space
$AP(\Pi, \mathbb R^n)\cap BC^1_t(\Pi,\R^n)$ into itself.

Let the condition $\bf{(B1)}$ be  fulfilled.
Due to the proof of   Theorem \ref{bounded_sol1}, this yields
$$
\bigl\|\widehat C(k)\bigr\|_{\LL\left(BC(\Pi;\R^n)\right)}+\bigl\|\widehat D(k)\bigr\|_{\LL(BC(\Pi;\R^n))}<1.
$$
  Hence, the operator $I -\widehat  C(k) - \widehat D(k)$ is bijective from $BC(\Pi;\R^n))$ into itself.
As a consequence, the solution 
$U^{k+1}\in BC(\Pi;\R^n)$ to the equation
$
U^{k+1}=(\widehat C(k) + \widehat D(k))U^{k+1}+\widehat F(k)(f,h)
$
is given by the (uniformly convergent) Neumann series
\beq\label{AP}
U^{k+1}=\left(I - \widehat C(k) - \widehat D(k)\right)^{-1}\widehat F(k)(f,h)=\sum\limits_{j=0}^\infty \left(\widehat C(k)+\widehat D(k)\right)^j\widehat F(k)(f,h).
\ee
Since the functions $f$ and $h$ are continuously differentiable in $t$, the function
$F(k)(f,h)$ belongs to $BC^1_t(\Pi,\R^n)$.
Moreover,
$$
\left(\widehat C(k)+\widehat D(k)\right)^j \mbox{ maps } AP(\Pi, \mathbb R^n)\cap BC^1_t(\Pi,\R^n) \mbox{ to }  AP(\Pi, \mathbb R^n)
$$
for each $j$.
Therefore, the  right hand side and, hence, the left hand side of \reff{AP}  belongs to 
$
AP(\Pi, \mathbb R^n)$. This means that that the function
$V^{k+1}=Q^kU^{k+1}$ belongs to $AP(\Pi, \mathbb R^n)$, as desired.

If the  condition $\bf{(B2)}$ (resp., $\bf{(B3)}$)
is fulfilled, then we use a similar argument. More precisely,
 we consider the formula \reff{uj} with $C$ and  $G_0$ (resp.,   with $C$ and  $H_0$)
replaced  appropriately by  $\widehat C(k)$ and  $\widehat G_0(k)$ (resp.,  by   $\widehat C(k)$ and  $\widehat H_0(k)$).
Taking into account  the inequalities
$\|\widehat G_0(k)\|_{\LL(BC(\Pi;\R^n))}<1$ (resp., $ \|\widehat H_0(k)\|_{\LL(BC(\Pi;\R^n))}<1
$),
we use
the Neumann series representation for the operator $(I-\widehat G_0(k))^{-1}$ (resp.,  $(I-\widehat H_0(k))^{-1}$)
to conclude  that the iterated solution $U^{k+1}$ and, hence 
$V^{k+1}=Q^kU^{k+1}$
belong to $AP(\Pi, \mathbb R^n)$. 
The proof is therefore complete.

\subsection{Proof of Part 2 of Theorem \ref{main}: Periodic solutions} \label{periodic}

We follow  the proof of the almost
periodic case in Section \ref{ap}, on each step referring to the periodicity instead of the almost periodicity. 
Obvious simplifications in the proof are caused by the identity
$\om_j^k(\eta,x,t) + T = \om_j^k(\eta,x,t+T)$.

\subsection{Proof of Theorem \ref{mainH}: Bounded solutions for space-periodic problems}\label{sec:boundedper}

The proof of Theorem \ref{mainH} repeats the proof of Theorem \ref{main}, with the only difference being that we need to refer to Theorem \ref{bounded_sol3} instead of Theorem \ref{bounded_sol1}.

\section*{Acknowledgments}
Irina Kmit and Viktor Tkachenko were supported by the VolkswagenStiftung Project ``From Modeling and Analysis to Approximation''.
Lutz Recke was supported by the DAAD program ``Ostpartnerschaften''.

\end{document}